\documentclass{article}

\usepackage{arxiv}

\usepackage[T1]{fontenc}    

\usepackage{url}            

\usepackage{nicefrac}       
\usepackage{microtype}      
\usepackage{amsmath}

\usepackage{graphicx}
\usepackage{indentfirst}
\usepackage{booktabs}
\usepackage{cite}

\usepackage[ruled,linesnumbered]{algorithm2e}
\usepackage{graphicx}
\usepackage[pagewise]{lineno}
\usepackage{xcolor,paralist,hyperref,fancyhdr,etoolbox}
\usepackage{tikz}
\usetikzlibrary{matrix, positioning, calc}
\usetikzlibrary{shapes, arrows.meta, positioning}
\usepackage[utf8]{inputenc}
\usepackage{csquotes}
\usepackage[english]{babel}
\usepackage{amsfonts}
\usepackage{multirow}
\usepackage{float}
\usepackage{cleveref}
\usepackage{algorithmic}
\usepackage[ruled,linesnumbered]{algorithm2e}
\newtheorem{theorem}{Theorem}

\makeatletter
\def\blx@maxline{77}
\makeatother

\hypersetup{ colorlinks=true, linkcolor=black, filecolor=black, urlcolor=black }

\usepackage{lipsum}

\title{A mechanism-driven reinforcement learning framework for shape optimization of airfoils}

  \author{
 Jingfeng Wang \\
  Department of Mathematics\\
  Faculty of Science and Technology\\
  University of Macau\\
  \texttt{shankswang953@gmail.mo} \\
   \And
 Guanghui Hu \\
  Department of Mathematics\\
  Faculty of Science and Technology\\
  University of Macau\\
  \texttt{garyhu@um.edu.mo} \\
}

\begin{document}
\maketitle
\begin{abstract}
  In this paper, a novel mechanism-driven reinforcement learning framework is proposed for airfoil shape optimization. To validate the framework, a reward function is designed and analyzed, from which the equivalence between the maximizing the cumulative reward and achieving the optimization objectives is guaranteed theoretically. To establish a quality exploration, and to obtain an accurate reward from the environment, an efficient solver for steady Euler equations is employed in the reinforcement learning method. The solver utilizes the Bézier curve to describe the shape of the airfoil, and a Newton-geometric multigrid method for the solution. In particular, a dual-weighted residual-based h-adaptive method is used for efficient calculation of target functional. To effectively streamline the airfoil shape during the deformation process, we introduce the Laplacian smoothing, and propose a Bézier fitting strategy, which not only remits mesh tangling but also guarantees a precise manipulation of the geometry. In addition, a neural network architecture is designed based on an attention mechanism to make the learning process more sensitive to the minor change of the airfoil geometry. Numerical experiments demonstrate that our framework can handle the optimization problem with hundreds of design variables. It is worth mentioning that, prior to this work, there are limited works combining such high-fidelity partial differential equatons framework with advanced reinforcement learning algorithms for design problems with such high dimensionality.
\end{abstract} 

\maketitle


\section{Introduction}
Aerodynamic design has been an important research theme for decades, particularly in the context of optimizing airfoils for enhanced aerodynamic performance\cite{martins2013multidisciplinary,nadarajah2007optimum} and fuel efficiency\cite{panagiotou2020aerodynamic}. 
Multidisciplinary analysis and optimization play an indispensable part in these fields.
Nonetheless, it is pointed out by an important proposal from NASA \cite{slotnick2014cfd} that high-fidelity computational fluid dynamics(CFD) is not routinely used in this area. Both academia and industry anticipate significant advancements in the integration of CFD methodologies by the year 2030. Aiming at conducting shape optimization of airfoils based on a high-fidelity CFD solver with the state-of-the-art technique, a mechanism-driven reinforcement learning framework is proposed in this paper.

 The traditional method pioneered by Jameson \cite{jameson1988aerodynamic, reuther1996aerodynamic} combined the optimization problem with adjoint equations which has been proven to be very efficient.
 However, the analysis and application of optimization techniques in the industry can not be generalized to more complicated target functional easily\cite{slotnick2014cfd}.
Combining CFD methods with optimization offers some distinct advantages. For instance, a more precise depiction of flow fields enables a more accurate calculation of target functionals. The advent of techniques such as the dual-weighted residual(DWR)-based mesh adaptation\cite{fidkowski2011review,DOLEJSI2021178,HARTMANN2015754} allows for the attainment of more precise target functionals while preserving computational accuracy. Moreover, the enhanced coupling between more accurate geometric representations and PDEs provides superior tools for simulating airfoil geometric deformations through the use of parametric curves and unstructured mesh methods to describe shapes\cite{meng2021fourth,Meng2022}.
 However, PDEs-constrained optimization also presents various challenges, including the high computational cost associated with solving PDEs\cite{biegler2003large} with complex geometries, non-linearities inherent in many physical systems, and application issues for ensuring that optimized shapes are physically realizable\cite{allaire2004structural}.

With the development of artificial intelligence in recent years, associating the neural networks with traditional computational fluid dynamics problems becomes feasible in shape optimization of airfoil\cite{li2022machine,chen2023study}. For example, the Bézier-GAN model in \cite{chen2020airfoil,chen2019aerodynamic} and the conditional variational autoencoder \cite{yonekura2021data} have been adopted to assess the performance of the shape of airfoils. Moreover, the recent advancements in reinforcement methodologies, such as deep reinforcement learning\cite{arulkumaran2017deep}, have significantly expanded its applicability and effectiveness in handling high-dimensional state and action spaces.

The distinctive advantages of reinforcement learning, compared with other machine learning algorithms, lie in its ability to make sequential decisions and learn optimal policies directly from interaction with the environment\cite{mnih2015human}. This is especially pertinent in the realm of control problems governed by partial differential equations\cite{duriez2017machine,rabault2019artificial}, where the dynamics are complex and highly nonlinear. Different from traditional machine learning approaches, reinforcement learning thrives in an exploration-based learning environment, making it uniquely suited for optimizing shapes and controls in PDEs-driven systems without the need for extensive pre-processed data.
 These developments, combined with our group's work in exploring advanced PDE solvers\cite{hu2010jcp,hu2011robust,hu2013adaptive,meng2021fourth}, have laid the groundwork for implementing a PDE solver-based reinforcement learning framework. 

Even though reinforcement learning has made significant strides in various applications, its integration into shape optimization, particularly within the context of PDEs-based systems, remains nascent. 
Ever since in \cite{lampton2010reinforcement}, reinforcement learning has been applied to morphing within an airfoil design context.  Initially, the application was constrained by the consideration of only four variables for airfoil morphing. With the advent of computational advancements, a direct shape optimization methodology utilizing reinforcement learning was introduced in \cite{viquerat2021direct}. Given that numerical simulations in CFD are time-consuming, the number of parameters is still far from the practical requirements. A transfer learning technique \cite{yan2019aerodynamic} is implemented to expedite the learning process. However, such a technique is still not a robust strategy without a mechanism involved.
In our research\cite{hu2016adjoint,HU2016235,wang2023towards,wang2024towards}, we have integrated the DWR-based mesh refinement strategy, a method that can generate a high quality mesh for the calculation of target functional. It saves time significantly, which supports the subsequent framework conducted with a large number of parameters.

Challenges in developing this framework combining reinforcement learning and PDE solvers are from the following aspects.

\textbullet ~The absence of efficient automation techniques in DWR-based mesh adaptation processes leads to critical challenges, particularly as significant geometric deformations occur.
    During the iterative process of reinforcement learning, the PDE solver might be invoked thousands of times. Manual intervention at each iteration is impractical\cite{lecun2015deep}. In \cite{nemec2014toward, wang2024towards}, the automation of DWR-based mesh adaptation is discussed. Further improvements should be implemented to guarantee a highly robust strategy even when the geometry is extraordinarily changed. Otherwise, the accuracy cannot be preserved, which will affect the agent's learning. 

\textbullet ~The deformation conducted within the unstructured mesh will bring potential mesh tangling issues\cite{johnson1994mesh}.
The change of geometry may be stochastic initially, where the quality of meshes can not be guaranteed. Besides, the curvature may change rapidly along the deformed shape. If the curvature variation occurs within the elements, it is challenging to capture the geometry correctly.

\textbullet ~Tackling the issues posed by controlling geometries with a large number of parameters is challenging\cite{chen2022manifold, perez2016automatic}. 
    The primary difficulty with a high parameter count is that it complicates the simulation process. Furthermore, while a more precise control over geometry enhances manipulability, it also poses challenges to maintaining a regular global shape.

\textbullet ~Designing a framework that efficiently meets varied optimization goals is inherently complex, requiring the creation of mechanisms specifically crafted for each unique task.
    The effectiveness of DWR-based mesh adaptation methods is inherently related to the current state of theoretical advancements\cite{fidkowski2011review}. The challenge is further compounded when dealing with more complex target functionals, whose resolution poses significant difficulties. Furthermore, translating traditional optimization problems into reinforcement learning frameworks introduces additional complexities. Ensuring the equivalence of primal optimization objectives and the result of reinforcement learning is paramount. This requires the analysis of the reward function\cite{dewancker2016bayesian}, necessitating it to be closely compatible with the underlying optimization context.

In this paper, we propose a novel mechanism-driven reinforcement learning, aiming at enhancing the shape optimization of airfoils. This approach effectively manages high-dimensional geometric parameters over unstructured meshes for different objectives, addressing several key challenges in the field. Prior to this work, there has been limited work combining the Twin Delayed Deep Deterministic Policy Gradient (TD3) algorithm with traditional PDE solvers for direct airfoil design optimization in such high-dimensional spaces. It is worth mentioning that we have analyzed the equivalence of the learning process and objectives based on specific reward function. 

We integrate an automatic algorithm leveraging the DWR-based mesh adaptation method, facilitating automation essential for iterative optimization and ensuring time efficiency without compromising accuracy.
Subsequently, we mitigate potential mesh tangling during airfoil deformation by incorporating Laplacian smoothing, which prevents irregular triangle accumulation and reduces finite element analysis errors due to mesh irregularities.
Then, the geometric manipulation utilizes up to 132 sample points for precise actions, coupled with a regularization strategy for control points generation. This approach ensures regular shape transitions and captures intricate curvature details, enhancing geometrical fidelity through Bézier fitting for effective dimensionality reduction.
Lastly, our framework employs the enhanced TD3 algorithm, tailored for high-dimensional control problems. 
It supports the optimization task with a well-aligned reward function, promoting behaviors yielding higher rewards and penalizing undesirable actions, thereby streamlining the learning process towards achieving optimization goals efficiently. 

The rest of this paper is organized as follows. In \cref{sec:basic}, we
briefly introduced the PDEs-constrained optimization, fundamental notations and the architecture of the
automatic DWR-based mesh adaptation method.  In \cref{sec:rl}, a basic process of the reinforcement learning based on the shape optimization of airfoil is introduced. Additionally, issues for setting up the learning process which originate from the high dimensionality and unstructured mesh are addressed. In \cref{sec:td3}, the algorithm of Twin Delayed DDPG is constructed. Further details of each component of the algorithm, analysis of the equivalence between the learning process and optimization task are elaborated.
In \cref{sec:result}, we present the numerical results. The future improvements are discussed in \cref{sec:conclusion}.

\section{PDEs-constrained optimization}
\label{sec:basic}
This study focuses on the shape optimization of airfoils using PDEs-constrained optimization techniques. Consider an objective function $J(\mathbf{u})$ that we seek to minimize, subject to a constraint imposed by a governing equation 
$R(\mathbf{u})=0$. The optimization problem can be mathematically formulated as follows:
\begin{equation}
\min_{\mathbf{u}} J(\mathbf{u}) \quad \text{subject to} \quad R(\mathbf{u}) = 0.
\end{equation}
In this work, we integrate PDEs-constrained optimization within a reinforcement learning framework. Our approach involves navigating the complex landscape of PDEs constraints to accurately achieve the desired target functional. We utilize the DWR method for mesh adaptation, which will be elaborated upon in this section. This method enhances our ability to refine the mesh dynamically based on error estimates that contribute directly towards improving the accuracy of the objective function $J$.

\subsection{Steady Euler equations}
In our previous works \cite{HU2016235, hu2016adjoint, hu2011robust}, we have developed a solver capable of handling the steady Euler equations. This solver not only facilitates the computation of these equations but also ensures the accurate resolution of the target functional, laying a solid foundation for our shape optimization endeavors. A brief introduction will be given in this section.

The steady Euler equations help us understand inviscid flow behavior and are capable of addressing the shape optimal design problem\cite{JamEuler}. These equations are represented in a conservative form for the inviscid two-dimensional steady state,
\begin{equation}
    \nabla\cdot\mathcal{F}(\mathbf{u}) = 0,\qquad \textbf{in } \Omega.
\end{equation}
Here $\mathbf{u}$ and $\mathcal{F}(\mathbf{u})$  are the conservative
variables and fluxes, which are given by
\begin{equation}
\mathbf{u} = \begin{bmatrix}
\rho \\ \rho u_x \\ \rho u_y \\ E
\end{bmatrix},
\quad \text{and} \quad \mathcal{F}(\mathbf{u}) = \begin{bmatrix}
\rho u_x & \rho u_y \\
\rho u_x^2 + p & \rho u_x u_y \\
\rho u_x u_y & \rho u_y^2 + p \\
u_x(E + p) & u_y(E + p)
\end{bmatrix},
\end{equation}
where the terms $(u_x, u_y)^T$, $\rho$, $p$, and $E$ represent the velocity, density, pressure, and total energy, respectively. The system is completed with the equation of state for an ideal gas, defined as
\begin{equation}
E = \frac{p}{\gamma - 1} + \frac{1}{2} \rho(u_x^2 + u_y^2),
\end{equation}
where $\gamma = 1.4$ indicates the specific heat ratio.

For instance, we consider the optimization of airfoil shapes with the primary objective of drag minimization. We define the drag as a specific quantity of interest \( J_d(\mathbf{u}) \), which is to be minimized. The mathematical formulation of the optimization problem is given by:
\begin{equation}
    \min_{\mathcal{S}} J_{d}(\mathbf{u}) \quad \text{subject to} \quad \nabla \cdot \mathcal{F}(\mathbf{u}) = 0, \quad \text{in} \ \Omega_{\mathcal{S}},
\end{equation}
where \(\mathcal{S}\) denotes the shape of the airfoil, and \(\Omega_{\mathcal{S}}\) represents the fluid domain bounded by the airfoil shape \(\mathcal{S}\). The governing equation \( \nabla \cdot \mathcal{F}(\mathbf{u}) = 0 \) encapsulates the flow dynamics within the domain, ensuring that the flow is modeled accurately according to the physical constraints.

Solving the PDEs and obtaining a precise target functional is our first concern. Then we discretize the domain $\Omega_{\mathcal{S}}$. Consider $\Omega_{\mathcal{S}}$ as a domain in $\mathbb{R}^2$ with boundary denoted by $\Gamma_{\mathcal{S}}$. We divide $\Omega_{\mathcal{S}}$ into distinct control volumes, $K_i$, via a shape-regular subdivision. Then the mesh under the current mesh size is denoted as $\mathcal{K}_h$. Each control volume is referred as $K_i$. The intersection between any two elements $K_i$ and $K_j$, is labeled as $e_{i,j}$, satisfying $e_{i,j} = \partial K_{i} \cap \partial K_{j}$. The edge $e_{i,j}$ relative to $K_{i}$ is equipped with a unit outer normal vector, $n_{i,j}$.
This setup allows us to redefine the weak form of Euler equations through the divergence theorem as:
\begin{equation}
    \mathcal{A}(\mathbf{u})=\int_{\Omega_{\mathcal{S}}} \nabla \cdot \mathcal{F}(\mathbf{u})dx=\sum_{i}\int_{K_i}\nabla\cdot
    \mathcal{F}(\mathbf{u})dx=\sum_{i}\sum_{j}\oint_{e_{i,j}\in\partial K_i}\mathcal{F}(\mathbf{u})\cdot n_{i,j}ds=0.
\end{equation}

To further advance towards a fully discretized system, we introduce a numerical flux function $\mathcal{H}(\cdot,\cdot,\cdot)$, given by:
\begin{equation}
    \sum_{i}\sum_{j}\oint_{e_{i,j}\in\partial K_i}\mathcal{H}(\mathbf{u}_i,\mathbf{u}_j, n_{i,j})ds=0.
    \label{EulerDiscrete}
\end{equation}

To tackle the nonlinear Equation \eqref{EulerDiscrete}, we employ the Newton-GMG framework. Details can refer to \cite{li2008multigrid, hu2010jcp,hu2013adaptive}. 
To obtain a precise value of the target functional, we turn to use the DWR-based mesh adaptation as part of the PDE-solver.

\subsection{An automatic DWR-based adaptation for the target functional calculation}
\label{sec:Auto}
In \cite{hu2016adjoint,wang2023towards}, we elaborated a framework innovated by Darmofal \cite{venditti2000adjoint,venditti2002grid}. To incorporate this technique within our framework, we present an overview of the methodology. 
Denote solutions on the coarse mesh as $\mathbf{u}_H$ and the finer mesh as $\mathbf{u}_h$. The initial goal is computing the integration of $J(\mathbf{u})$ over a fine mesh, denoted as $J_h(\mathbf{u}_h)$. A multi-variable Taylor series expansion method is adopted as follows:
\begin{equation}
  J_h(\mathbf{u}_h) = J_h(\mathbf{u}_h^H) + \left.\frac{\partial J_h}{\partial \mathbf{u}_h}\right|_{\mathbf{u}_h^H}(\mathbf{u}_h-\mathbf{u}_h^H) + \cdots,
\end{equation}
where $\mathbf{u}_h^H$ denotes the coarse mesh solution $\mathbf{u}_H$ projected onto the fine mesh space $\mathcal{V}_h$ through a prolongation operator $I_h^H$, $\mathbf{u}_h^H = I_h^H\mathbf{u}_H.$

We denote the residual of governing equations in the space $\mathcal{V}_h$ by $\mathcal{R}_h(\cdot)$, which satisfies the equations $\mathcal{R}_h(\mathbf{u}_h) = 0.$ Upon linearizing, we derive
\begin{equation}
  \mathcal{R}_h(\mathbf{u}_h) = \mathcal{R}_h(\mathbf{u}_h^H) + \left.\frac{\partial \mathcal{R}_h}{\partial \mathbf{u}_h}\right|_{\mathbf{u}_h^H}(\mathbf{u}_h-\mathbf{u}_h^H) + \cdots.
\end{equation}
Given $\left.({\partial \mathcal{R}_h}/{\partial \mathbf{u}_h})\right|_{\mathbf{u}_h^H}$ represents the Jacobian matrix, its inversion approximates the error vector as,
\begin{equation}
  \mathbf{u}_h - \mathbf{u}_h^H \approx -\left.({\partial \mathcal{R}_h}/{\partial \mathbf{u}_h})^{-1}\right|_{\mathbf{u}_h^H}\mathcal{R}_h(\mathbf{u}_h^H).
\end{equation}
Substituting the error vector into the residual equation yields the updated quantity of interest,
\begin{equation}
  J_h(\mathbf{u}_h) = J_h(\mathbf{u}_h^H) - (\mathbf{z}_h|_{\mathbf{u}_h^H})^T\mathcal{R}_h(\mathbf{u}_h^H),
\end{equation}
where $\mathbf{z}_h|_{\mathbf{u}_h^H}$ is derived from the fully discrete dual equations:
\begin{equation}
  \left(\left.\frac{\partial \mathcal{R}_h}{\partial \mathbf{u}_h}\right|_{\mathbf{u}_h^H}\right)^T\mathbf{z}_h|_{\mathbf{u}_h^H} = \left(\left.\frac{\partial J_h}{\partial \mathbf{u}_h}\right|_{\mathbf{u}_h^H}\right)^T.
\end{equation}

The calculation of $\mathbf{z}_h$ on the embedded mesh was expensive. The literature presents various methods to streamline this process. For instance, the works of \cite{venditti2000adjoint, nemec2007adjoint} introduce an error correction method that interpolates dual solutions onto a fine mesh, while the BiCG method in \cite{DOLEJSI2021178,dolejvsi2023anisotropic} has been adopted to calculate the dual and primal problems at once. With the help of machine learning techniques, we developed the hybrid CNNs-Dual framework\cite{wang2024towards} to accelerate the calculation of dual solutions which is very effective. However, the CNNs-Dual solver can not preserve the accuracy for deformed geometry. 
 Even though techniques like GCN(Graph Convolutional Networks)-enhanced learning can mend this issue, 
 we still adopt the DWR-based adaptation developed in \cite{wang2023towards} in this work for its robustness.

\begin{figure}[!h]\centering
\frame{\includegraphics[width=0.48\textwidth]{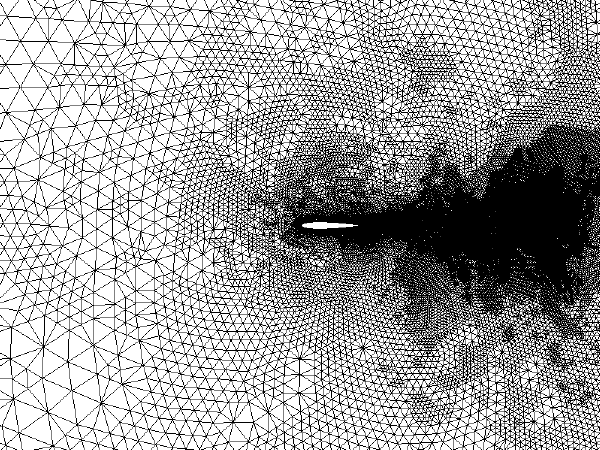}}
        \frame{\includegraphics[width=0.48\textwidth]{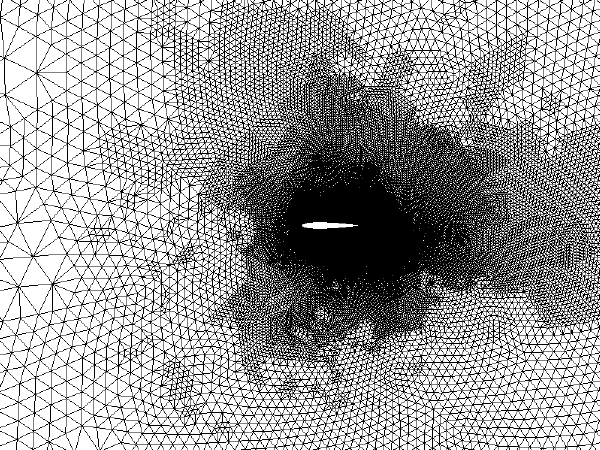}}
         
       \caption{Mesh around the NACA0012 airfoil. Left: the residual-based mesh adaptation. Right: the DWR-based mesh adaptation.}
      \label{DifferentAdapt}
      \end{figure}

Although the dual solutions can be solved effectively, issues persist in balancing the value of primal residuals and dual solutions as their magnitudes are not compatible. Indicators for each element are generated based on their multiplication $
      \eta_{K_i} = \sum_{K_{ij}}|(\mathbf{z}_h)^T\mathcal{R}_h(\mathbf{u}_h^H)|,
$
    where $K_{ij}$ is the subelements $K_{i}$.
    For the implementation of mesh adaptation, a parameter needs to be set for guiding such a process. However, setting this parameter is tricky because a lower value may lead to an approximately uniform refinement while a higher value may not solve the target functional with expected accuracy. What makes it worse is that such a parameter differs with the update of configurations and different geometries. Thus, it is crucial to find an automatic parameter-setting strategy for guaranteeing the mesh adaptation process during the optimization. In \cite{wang2024towards}, we developed such a technique based on statistical analysis. At last, the procedure of the algorithm for calculating the target functional is concisely organised as  \cref{originalAlgorithm}:

\begin{algorithm}
\caption{DWR-based adaptation for calculating target-functional}
\label{originalAlgorithm}
\begin{algorithmic}[1]
\REQUIRE Initial mesh $\mathcal{K}_H$
\ENSURE Refined mesh $\mathcal{K}_h$ and target functional $J_h(\mathbf{u}_h)$
\WHILE{Iter $<$ refine\_step}
    \STATE Set the current mesh as $\mathcal{K}_H$
    \STATE Solve $\mathcal{R}_{H}(\mathbf{u}_H)=0$ with a residual tolerance of $1.0\times 10^{-3}$
    \STATE Interpolate solution $\mathbf{u}_H$ from mesh $\mathcal{K}_{H}$ to $\mathcal{K}_{h}$ to obtain $\mathbf{u}_h^H$
    \STATE Compute the residual $\mathcal{R}_h\left(\mathbf{u}_h^H\right)$
    \STATE Solve $\mathcal{R}_{h}(\mathbf{u}_h)=0$ with a residual tolerance of $1.0\times 10^{-3}$
    \STATE Generate the dual solutions from the primal residuals
    \STATE Obtain $\mathbf{z}_h$ by solving the dual equation with GMG
    \STATE Determine the error indicator for each element
    \STATE Automatic strategy for selecting $TOL$
    \WHILE{any element $K_H$ has $\mathcal{E}_{K_H}>TOL$}
        \STATE Refine element $K_{H}$ adaptively as described in \cite{HU2016235}
    \ENDWHILE
\ENDWHILE
\STATE Obtain the target functional $J_h(\mathbf{u}_h)$
\end{algorithmic}
\end{algorithm}

This efficient algorithm ensures the accuracy of the target functional can be maintained across configurations and geometries not typically encountered in engineering problems. During the exploration process inherent in the reinforcement learning framework, scenarios may arise that deviate significantly from those within the training datasets. Under these circumstances, a learning-based solver generally struggles to maintain accuracy, potentially feeding incorrect experiences back to the agent. In contrast, the paradigm of PDEs-constrained optimizing via reinforcement offers a more precise framework to assist the agent in making judgments.

\section{A reinforcement learning framework for shape optimization of airfoils}
\label{sec:rl}
The iteration process of reinforcement learning can be treated as an agent interacting with the environment and accumulating experience. The shape of the airfoil, composed of hundreds of points, can be conceptualized as a vector in a high-dimensional space, essentially representing a point in this space. Thus, all potential airfoil shapes can form a set within this high-dimensional space. The objective function, in this context, is a functional defined on this set. The goal of the reinforcement learning process is to locate the extremum on the set. To discuss the learning framework for the shape optimization of airfoil, the Markov decision process(MDP) with specific components should be defined correspondingly. However, the mesh deformation conducted within the unstructured mesh is challenging, where mesh tangling, curvature oscillation and irregular shapes are potential issues that will greatly influence the learning process. In this section, we will introduce a basic process of the reinforcement learning and address these potential issues.

\subsection{Markov decision process}
\begin{figure}[!h]\centering
        \includegraphics[width=0.8\textwidth,height=0.3\textheight]{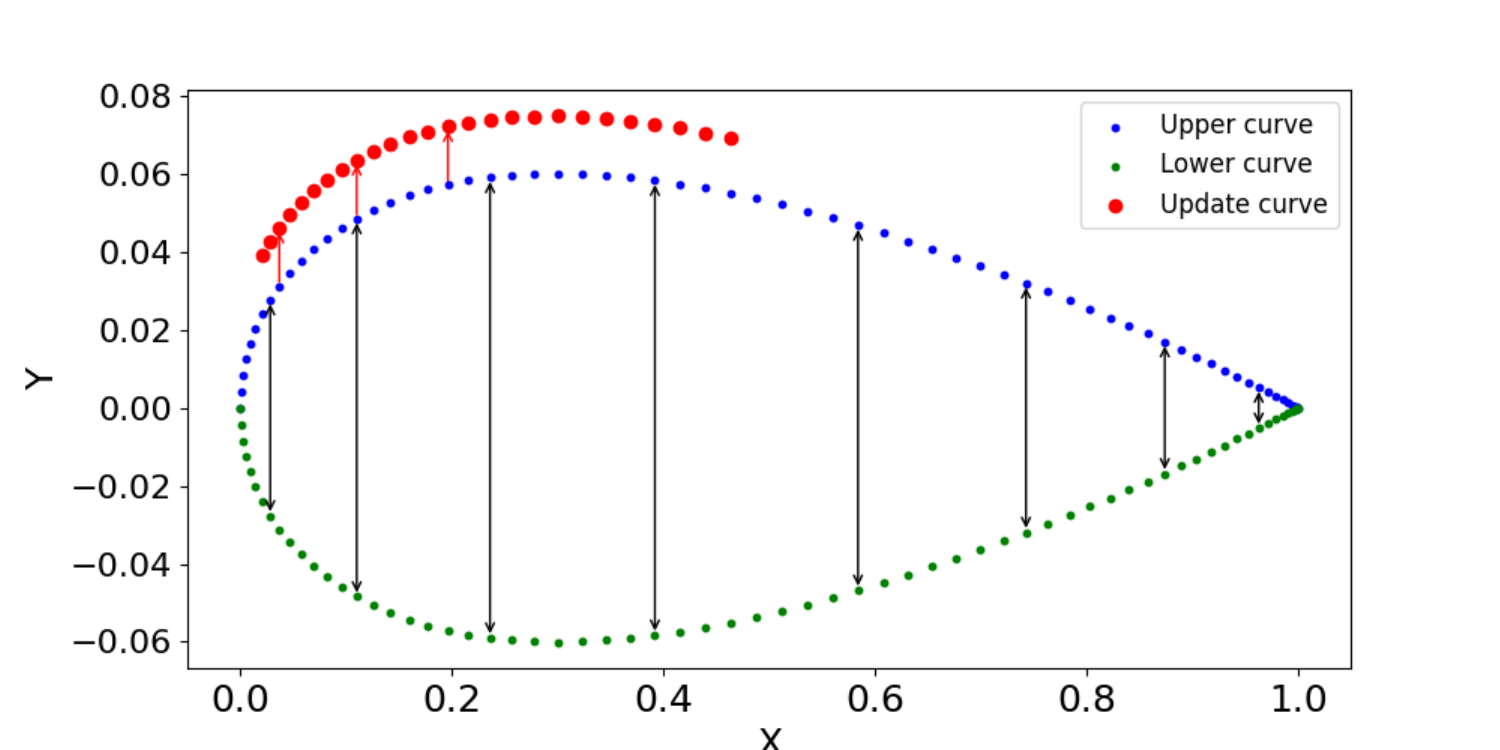}
       \caption{Update action with a smooth function acted on the geometry with $132$ sample points.}
      \label{smoothAction}
      \end{figure}

The Markov decision process can be defined as follows.
 At each time step $t$, a 4-tuple $(S_t, A_t, P_t, R_t)$ is adopted to record the process \cite{bellman1957markovian} which contains the interaction between agent and environment, where $t = 0, 1, 2, \cdots.$ Here $S_t$ represents the state, consisting of the mesh points that depict the geometry of the airfoil at time step $t$. $A_t$ stands for the action performed on $S_t$. The action $A_t$ consists of three parameters $[x_{target}, y_{Uchange}, y_{Lchange}]$. Given $x$ coordinate $x_{target}$, there exists a pair of points on both curves whose values on $y$ direction are $y_{upper}$ and $y_{lower}$. Then the $y$ coordinates will be updated with parameters $ y_{Uchange}$ and $y_{Lchange}$ on both curves respectively. However, the values will not only be acted with respect to $x_{target}$. This behavior introduces a sharp deformation, which will bring irregularity to the geometry and may even cause a mesh tangling phenomenon. To perform a more smooth and natural action on the airfoil, we utilize a smooth technique as shown in Figure \ref{smoothAction}. There are totally 132 sample points to depict the shape of the airfoil. 
 For a specific $x_{target}$, the update value $y_{Uchange}$ will not be acted on the point directly. 
The delta function, defined as 
$
    \displaystyle SA(x) = e^{-(x-x_{target})^2/2\delta^2},
$
will be set as a coefficient that can act on the neighbor sample points and smoothen the deformation process. $\delta$ is set as $0.2$ to $0.8$ in this model. Then the $y$ coordinate of the upper curve will be updated with 
  $\displaystyle  y' = y + y_{Uchange}SA(x).$
  
After the deformation, our algorithm will generate the Bézier curves. The update points represent the new state of the airfoil, $S_{t+1}$. At each time step, there are various actions that the agent can take. $p(S_{t+1}| S_{t})$ is denoted as the transition probability for the state $S_t$ change to $S_{t+1}$ with action $A_t$. The agent will evaluate the potential reward from different actions. The action that the agent takes at state $S_{t}$ is $A_t$, which is generated from the policy $\pi_t : S_{t}\rightarrow A_t$. The agent will evaluate the reward $R_{t+1}(S_t, A_t, S_{t+1})$ from the environment and try to maximize the cumulative reward from a potential trajectory $\tau := (S_0, A_0, S_1, A_1, \cdots , S_t, A_t, \cdots).
$
In this trajectory, the cumulative reward at time step $t$ is defined as $
    G_t(\tau) = \sum_{i=0}^{T}\gamma^{i}R_{t+i+1}.
$
Then the agent shall be trained to maximize the cumulative reward so that the deformation can fulfill the optimization target. The Twin Delayed DDPG algorithm\cite{fujimoto2018addressing} is constructed in this work for its efficiency in continuous space and is suitable for the current framework. 
In the following part, we start with two important components during the learning process to understand such a learning task, actor and critic.

 \subsection{An actor-critic scheme} 
  The actor, parameterized by $\theta^a$, is responsible for determining the best possible action $A_t$, given the current state $S_t$ at time step $t$ with the current policy, i.e., $A_t = \pi_{\theta^a}(S_t)$. This direct policy optimization allows the agent to efficiently explore the continuous action space and identify optimal actions.

The actor will generate potential actions that form a state-action sequence, i.e., $\tau_{\theta^a} \sim \pi_{\theta^a}$. The cumulative reward at time step $t$ cannot explicitly calculated as the reward should be obtained from the CFD solver, which is generally time-consuming while the infinite sequence cannot be obtained. The critic, parameterized by $\theta^c$, evaluates cumulative reward based on corresponding states. Hence, $\tilde{R}_{t+1}(S_t, A_t| \theta^c)$ is the cumulative reward calculated by the critic with a neural network. For a specific long-term state evolution $\tau$, 
\begin{equation}
    \tilde{G}_t(\tau_{\theta^a}|\theta^c) = \tilde{R}_{t+1}(S_t, A_t| \theta^c) + \sum_{i=1}^{T}\gamma^{i}\tilde{R}_{t+i+1}\Big((S_{t+i} | \theta^a), (A_{t+i} | \theta^a)|\theta^c \Big)
\end{equation}
is defined as calculating the expected cumulative reward. The Q-value, which is represented as 
\begin{equation}
    Q(S_t|\theta^a, \theta^c) = \mathbb{E}_{\tau_{\theta^a} \sim \pi_{\theta^a}}\Big(\tilde{G}_t(\tau_{\theta^a}|\theta^c)\Big),
\end{equation}
 is adopted to assess the performance of the actor based on $\theta^a$.
 The critic guides the actor towards actions that maximize the long-term reward. The feedback of critic is essential for updating the parameters of actor and refining the policy.

 During training, a replay buffer is utilized to store and randomly sample previous transition tuples $(S_t, A_t, S_{t+1}, R_{t+1})$ , facilitating the learning of uncorrelated and diversified experiences. 
In the training process, the critic is updated by minimizing the mean squared error between a batch of current Q-values and target Q-values. The actor is then updated by performing gradient ascent on the expected Q-values output by the critic, thus steering the actor towards actions that yield higher rewards.

Through the interaction between the actor and the critic, the agent effectively navigates the continuous action space with the experience in the replay buffer. As a result, it can learn sophisticated policies and perform expected actions on the airfoil, thereby optimizing its performance in complex environments. However, performing actions on this problem, based on such a high-dimensional parameters within the unstructured mesh, is not straightforward. Further improvement should be conducted. 

\subsection{Action of mesh deformation on unstructured mesh}
\begin{figure}[h]\centering
        \includegraphics[width=1.0\textwidth,height=0.35\textheight]{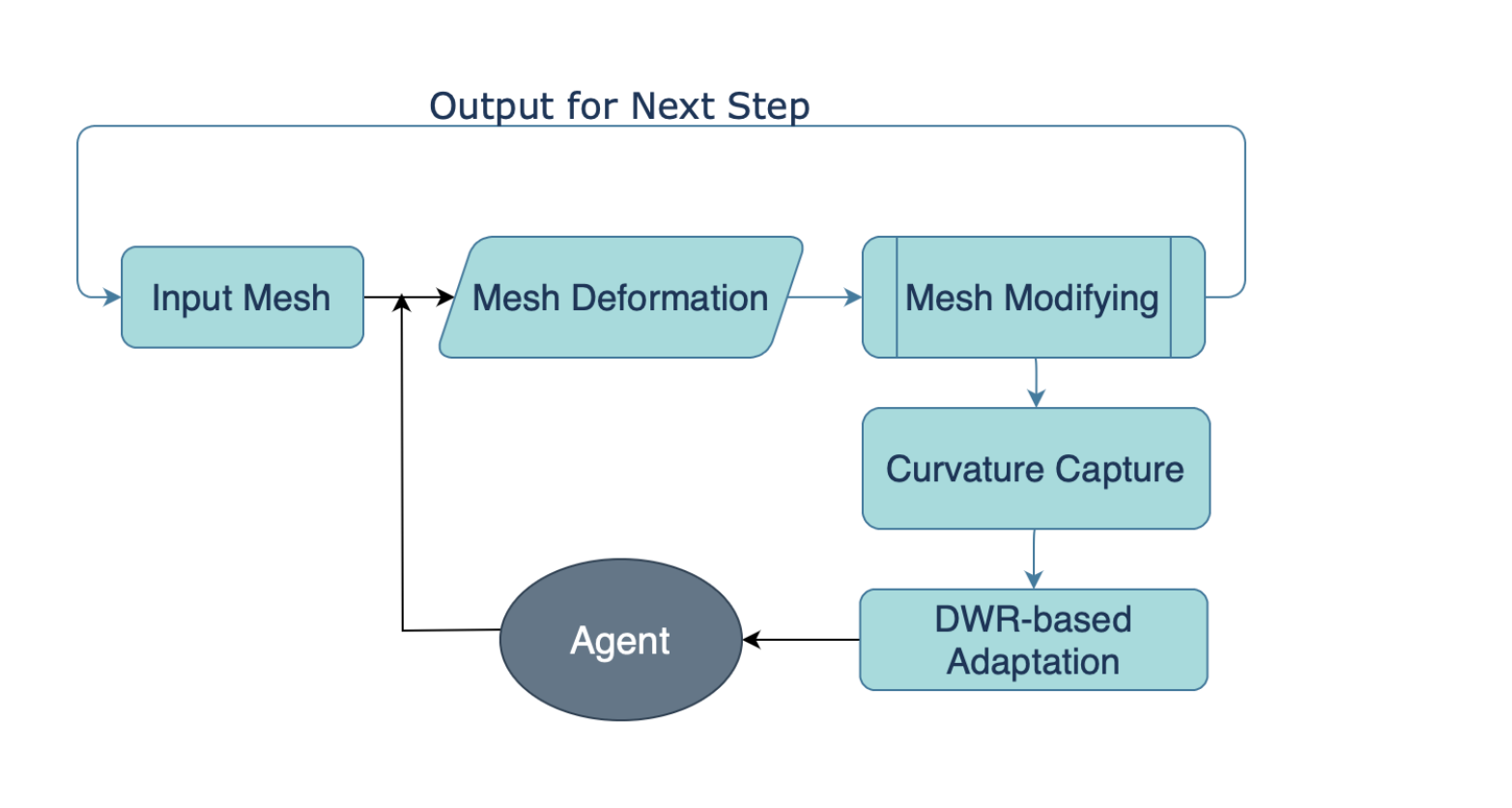}
      \caption{Mesh deformation, modification, and calculation process during the reinforcement learning}  
      \label{RL4deform}
      \end{figure}

Mesh deformation during shape optimization is a complicated process. It is important to ensure the quality of the mesh is capable of calculating target functional. The process of mesh deformation, modification, and calculation can be organized as  Figure \ref{RL4deform}.

To conduct the deformation of the airfoil, we need to adopt the parameterized curve for depicting the geometry at first. In this study, we applied the Bézier curves representation for the parameterization. The airfoil consists of an upper curve and a lower curve. For a given $x$ coordinate, there exists a pair of $y$ coordinates on the two parts. The representation for such two parts can ensure the $y$ coordinate of the upper half is always larger than that of the lower part so that the two curves will not intersect and prevent the mesh from tangling.

With the Bézier least square fitting, control points can be calculated with the sample points along the airfoil. Then the Bézier curves can be generated by the Bernstein polynomials, i.e.,
\begin{equation}
    B(t) = \sum_{i=0}^{n} \binom{n}{i} (1-t)^{n-i} t^i CP_i.
\label{BernsteinPolynomial}
\end{equation}
Here $B(t)$ represents the Bézier curve, which is a function of the parameter $t$. $t$ ranges from $0$ to $1$ which essentially interpolates between the control points $CP_i$. Then the airfoil with the defined curve is set in a circle with a large enough radius. With such a technique, we defined the computational domain on an unstructured mesh. The airfoil is depicted with a determinate number of mesh points. However, difficulties will brought by the unstructured mesh and high dimensionality with the iteration get processed. We will elaborate these issues as follows.

\begin{figure}[!h]\centering
        \frame{\includegraphics[width=0.48\textwidth]{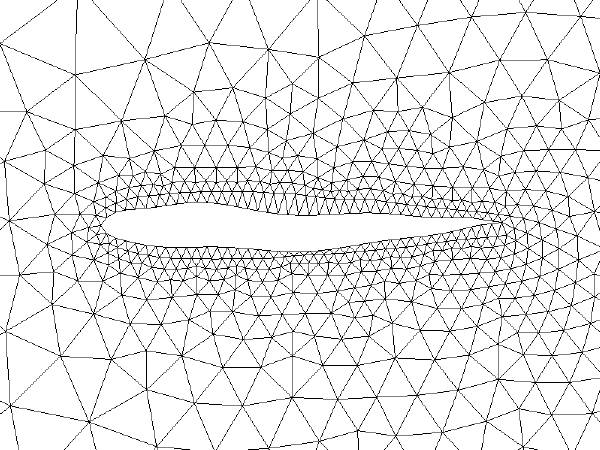}}
        \frame{\includegraphics[width=0.48\textwidth]{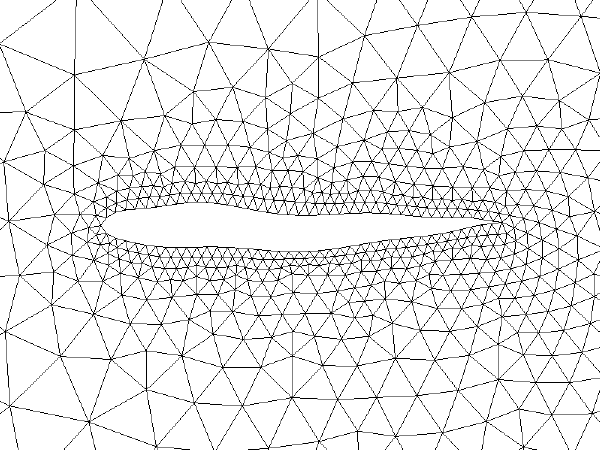}}
       \caption{Meshes around the airfoils. Left: before  modifying; Right: after modifying; }
      \label{Meshmodify}
      \end{figure}

\textbullet ~ \textbf{Mesh tangling.}
When the agent tries to modify the shape of airfoils, it will change all the sample points. Consequently, the control points will be updated and result in varied Bézier curves. From the updated expression of Bézier curves, we can calculate the mesh points along the airfoils in an anticlockwise direction. Afterward, we build a map structure so that the mesh file can find the corresponding mesh points on the airfoils and update their coordinates. As shown in the left part of Figure \ref{Meshmodify}, the geometry received the action from the agent and the shape gets updated. 

However, it is shown in the left part of Figure \ref{Meshmodify} that the quality of mesh gets influenced since the movement of boundary mesh points may lead to a sharp acute angle for the triangles near the upper curve and a large obtuse angle for the triangles near the lower half. While the optimization lasts for a long period, the elements around the airfoil will become more irregular with iteration steps, leading to the mesh tangling issues. Hence we shall modify the mesh so that it is capable of dealing multi-step mesh deformation.
The movement of points is restricted to a specific scale. Then the deformed mesh should be modified for further calculation. Since the most closely affected elements are those on the boundary part, such an area should be processed specially. For those boundary elements that have two points on the airfoil, the third point is relocated so that this element will be a regular triangle. Then we adopted the Laplacian smoothing technique to amend the quality of the mesh. The fundamental operation in Laplacian smoothing involves updating the position of a node to be the average of the positions of its neighboring nodes. Mathematically, this can be expressed as $\mathbf{x}_i^{'} = \frac{1}{m}\sum_{j=1}^{m} \mathbf{x}_j$,    
where $\mathbf{x}_i$ is the coordinate vector of the $i$-th point and $\mathbf{x}_i^{'}$ is the update value. $\mathbf{x}_j$ with index from $1$ to $m$ are neighboring nodes of the original point. With such a modification process, the quality of the mesh is enhanced. As shown in the right part of Figure \ref{Meshmodify}, the elements around the airfoil are recovered so that sharp acute angle and large obtuse angle can be prevented. 

Aiming to maintain efficiency throughout the iterative process, re-meshing is not considered as it is time-consuming, especially when generalized to three-dimensional cases. Even though the post-smoothing mesh possesses improved structural integrity of the triangular elements, it should not be far away from the pre-processed state. To this end, we leverage the C++ Standard Template Library's map data structure for its capability to manage modifications without drastically altering point positions or changing the number of vertices within the mesh.

In our approach, vertices that are not on the boundary but belong to boundary triangles are assigned with level 1. Subsequently, vertices without a level, belonging to triangles that include level 1 vertices, are marked as level 2. Then we recalculate the triangles containing level 1 and level 2 vertices so that the new triangles can approximate an equilateral triangle as possible. This recalculated position is mapped from the initial position, facilitating updates to the mesh data file through these mappings. Thus, without a re-meshing process, it allows for an effective improvement in mesh quality, crucial for the accuracy of subsequent simulations and analyses.

\textbullet ~\textbf{Curvature oscillation.}
It is worth noting that when the agent tries to modify the geometry, the curvature of the airfoil changes rapidly. The deformation conducted with a scale comparable to the mesh size can be handled with the method mentioned above. However, when the scale is smaller than the mesh size, the given mesh points can not detect such curvature variation within the segment. Then, it will add noise to both the primal and dual equations, and return an imprecise reward back to the agent.
As a result, it will bring more complexities and increase training difficulties. Out of this reason, we shall use triangles with a smaller scale. 

\begin{figure}[h]\centering
        \frame{\includegraphics[width=0.48\textwidth]{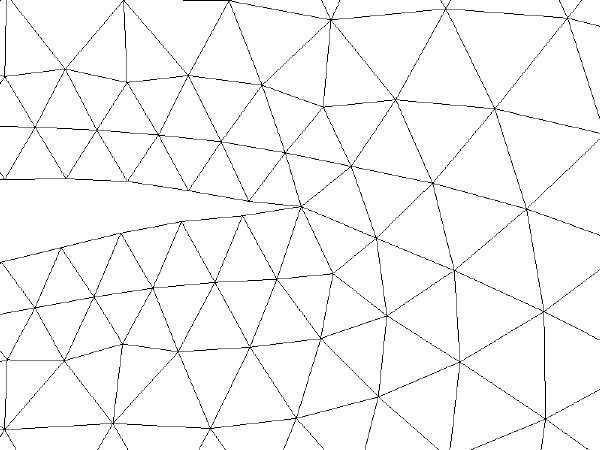}}
        \frame{\includegraphics[width=0.48\textwidth]{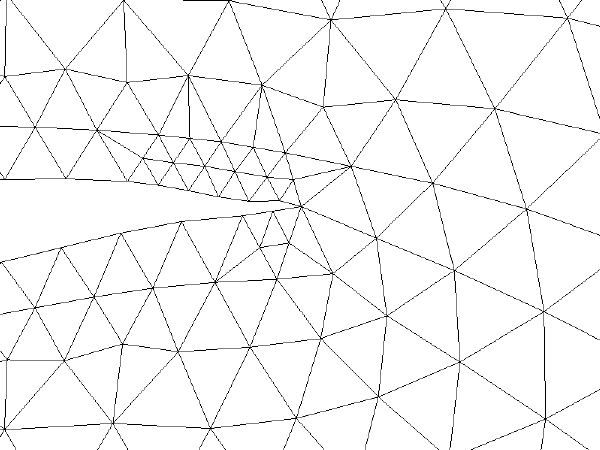}}
       \caption{Mesh around the tail. Left: before curvature capture; Right: after curvature capture.}
      \label{curvatureCap}
      \end{figure}
      
To calculate the curvature for each boundary curve, we shall find the parameters correspond with the endpoints, $t_1 = B^{-1}(\mathbf{x}_1)$ and $t_2 = B^{-1}(\mathbf{x}_2)$. Then we generate sample points that are dense enough within such segment, $sp_0, sp_1, \dots, sp_{L+1}$. For three consecutive sample points $sp_{i-1}, sp_{i}, sp_{i+1}$, we use $\theta_i$ to denote the angles between $\overrightarrow{sp_{i-1}sp_i}$ and $\overrightarrow{sp_{i}sp_{i+1}}$, then $\kappa_i ={\theta_i}\slash{\| \overrightarrow{sp_{i}sp_{i+1}} \|}$ is used to denote the discrete curvature values. $\sum_{i=0}^{L} \kappa_i$ is used to denote the curvature change within the segment.
It is shown in the left part of Figure \ref{curvatureCap}, that if the curvature change scale is smaller than the mesh size like the elements near the tail, the mesh can not detect such variation. However, when we use elements with smaller sizes, it can capture the curvature change within the segment. After running this algorithm several times, the mesh can depict a more elaborate airfoil shape. Even though issues remains on the irregular shapes as we try to manipulate the airfoil with hundreds of sample points.

\textbullet ~\textbf{Irregular shape.}
 High-dimensional state representations, while offering detailed control over the geometry, exacerbate the curse of dimensionality. This not only makes the optimization landscape more complex and harder to navigate but also strains the underlying PDE solvers. For machine learning-based solvers, deviations from the trained geometries can result in diminished accuracy, adversely affecting the fidelity of the objective functional evaluations. This necessitates the integration of traditional numerical methods capable of handling a wider variety of shapes with robustness, ensuring high precision in objective functional computation even when the configuration changed significantly.

\begin{figure}[!h]\centering
        \frame{\includegraphics[width=0.48\textwidth]{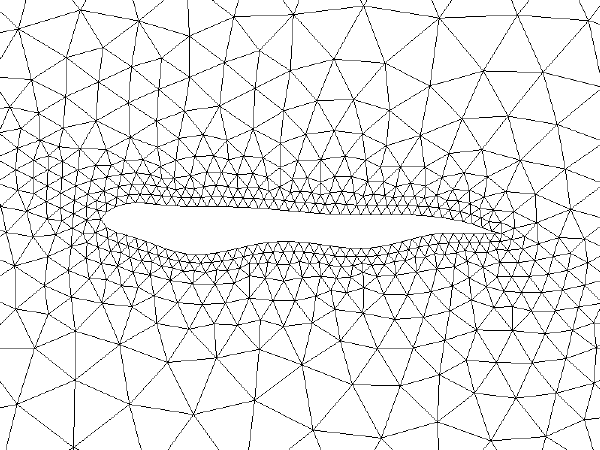}}
        \frame{\includegraphics[width=0.48\textwidth]{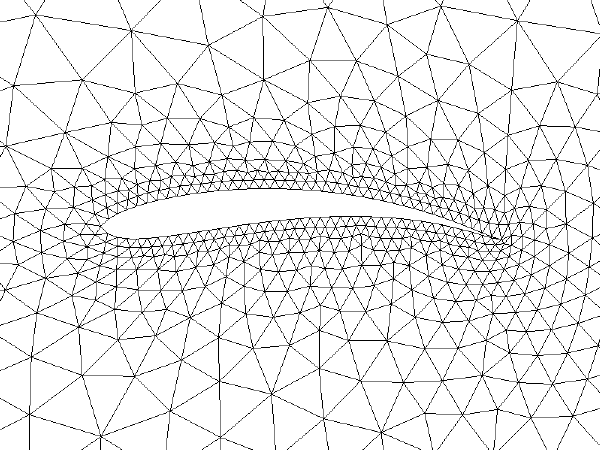}}
       \caption{Meshes around the airfoil. Left: without curve smoothing; Right: with curve smoothing; }
      \label{Meshsmoothing}
      \end{figure}

Furthermore, hundreds of sample points can introduce local irregularity to the geometry such as shown in Figure \ref{Meshsmoothing}. While actual engineering solutions rarely embody such erratic geometric characteristics, their presence during the learning process can significantly affect the optimization. Irregular geometries not only challenge the solver's capacity to provide accurate solutions but also obscure the learning signal, making it difficult for the agent to discern productive learning directions. Such geometric irregularities could mislead the agent, leading to suboptimal policies that fail to converge to the desired outcome. 

An improved approach entails adopting dimensionality reduction to extract essential geometric features into a more computationally efficient representation. While directly reducing the dimensionality of sample points will compromise their utility for precise geometric manipulation, strategies should be designed to take a balance between the local controllability and regularity. The utilization of the Bézier curve fitting stands out as a proficient technique for reducing dimensionality while capturing critical curve characteristics. This method effectively summarizes vital information, facilitating a streamlined representation of curves. However, the challenges still persist, particularly regarding the regularity of the curves. Directly lowering the order of Bézier curves can inadvertently eliminate valuable local details. To circumvent these issues, incorporating smoothing techniques during the Bézier curve generation process emerges as a potent solution. Integrating regularization directly into the curve fitting algorithm ensures the regularity of the curves, preserving the integrity of local geometric features without sacrificing the overall quality of the representation. 

The Bézier curve fitting for above purpose can be formulated by :
\begin{equation}
\begin{array}{cc}
     \text{Find}\quad  \{CP_i\}_{i=0}^{n}, \qquad
     \min \quad \sum_{j=0}^{L+1}|\!| B(t_j) - sp_j|\!|^2 ,
         \end{array}
\end{equation}
where $t_j$ is the parameters corresponding to the sample points $sp_j$. With a regularization term added to the system, the equation becomes:
\begin{equation}
\begin{array}{cc}
     &\text{Find}\qquad  \{CP_i\}_{i=0}^{n}, \\
     &\min \qquad \sum_{j=0}^{L+1}|\!| B(t_j) - sp_j|\!|^2 + \lambda_{s} * \sum_{i=1}^{n-1}|\!| CP_{i-1} - 2 * CP_i + CP_{i+1}|\!|^2
    \end{array},
\end{equation}
where $\lambda_s$ stands for the smoothing coefficient.
Such an regularization term is to make the control points equally distributed. Then the potential issues of the irregular curve can be mitigated effectively. As shown in  \ref{Meshsmoothing}, the deformed curve preserves a satisfactory geometry while having a detailed control of the airfoil.

\section{Design of Twin Delayed DDPG algorithm}
\label{sec:td3}
In the previous section, we introduced the learning process that can manipulate the deformation of the airfoil shape. However, training an agent to achieve the goal of interacting with the environment in such a complicated process is not straightforward. 
With the development of algorithm in reinforcement learning, the Twin Delayed DDPG framework enhances the training performance effectively.
In this section, we will elaborate on the details of the Twin Delayed DDPG framework for the design of airfoil.

\subsection{Actor and critic neural networks}

    \begin{figure}[h]\centering
        \includegraphics[width=1.0\textwidth,height=0.42\textheight]{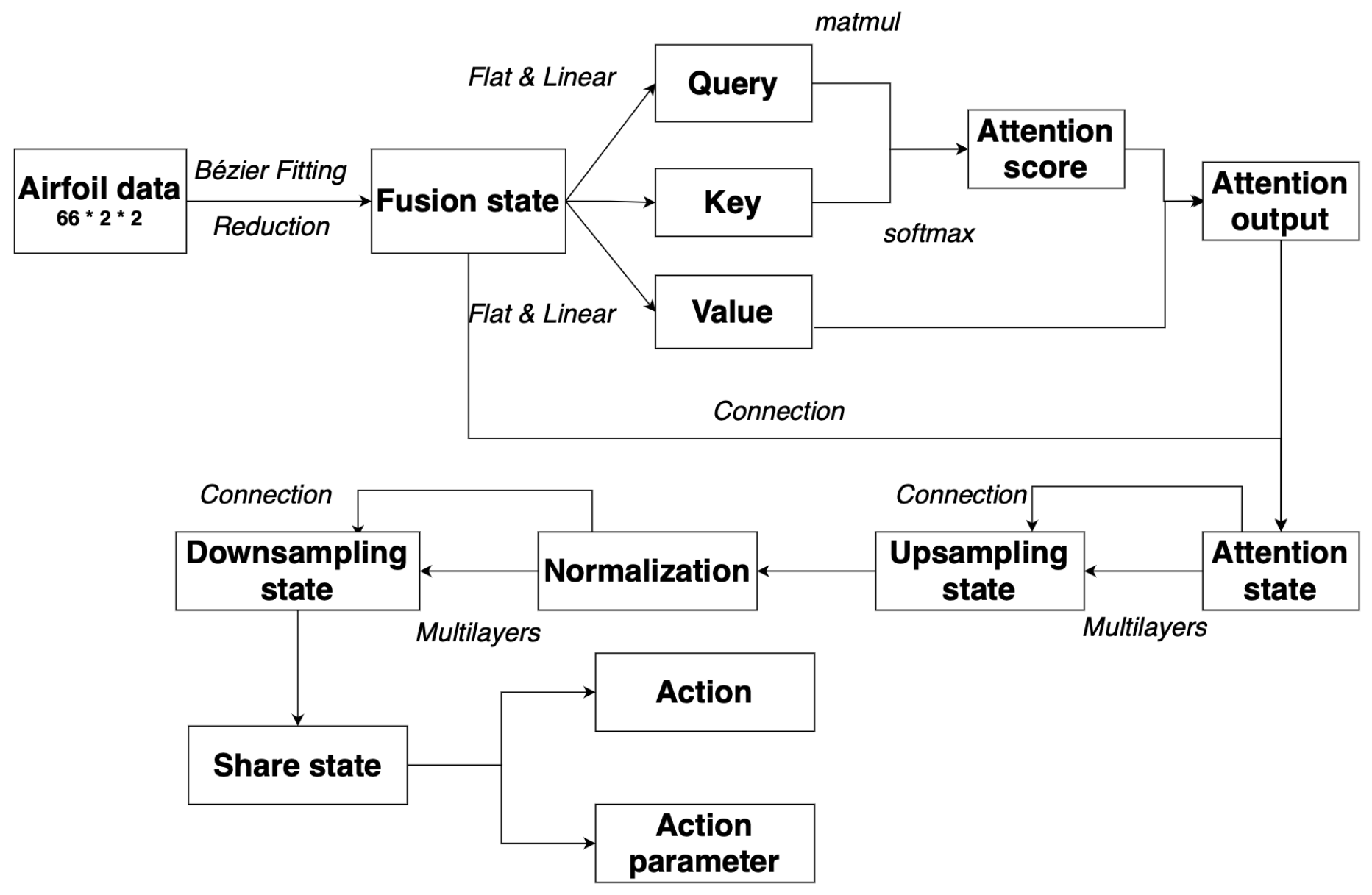}
      \caption{Actor network architecture}  
      \label{tab:actor_network_architecture}
      \end{figure}
      
The actor, as shown in \cref{tab:actor_network_architecture}, is developed to choose the potential action based on the input state such that the obtained reward is maximal. While the state is the tensor representation of the corresponding airfoil shape, the actor should capture the main feature and distinguish different shapes. Aiming at having a precise control of the airfoil, the geometry is constructed with hundreds of points pair. However, the agent may fail to learn valuable actions based on different states when this tensor representation is conducted within a high-dimensional space. A new state conducted by a specific action only takes a small change of target location. Such a difference is easy to overlook. The smooth action process explained in the previous section can remit such influence but the essential difficulties have not been solved. Trying to conduct the deformation within a lower dimensional tensor space, reduction module should be constructed.




 The variational autoencoder meets the requirement of dimensionality reduction and we tested the integration at first. However, it brings extra problems to the construction of the learning process. Firstly, it is not easy to take the balance of the reconstruction error and original loss function. Secondly, deeper neural networks aggravate the risk of the gradient disappearance problem. Thirdly, the sampling process is expensive in this reinforcement learning task. The VAE module can not perform well with limited data and gets trapped into overfitting issues. 

 Then we tried a different method, using Bézier curves representation for the dimensionality reduction. The Bézier curves take the most important information for specific curves. Besides, the difference between different control points originating from corresponding states is amplified with the least square process. It helps the agent to distinguish different states. 

    \begin{figure}[h]\centering
        \includegraphics[width=1.0\textwidth,height=0.4\textheight]{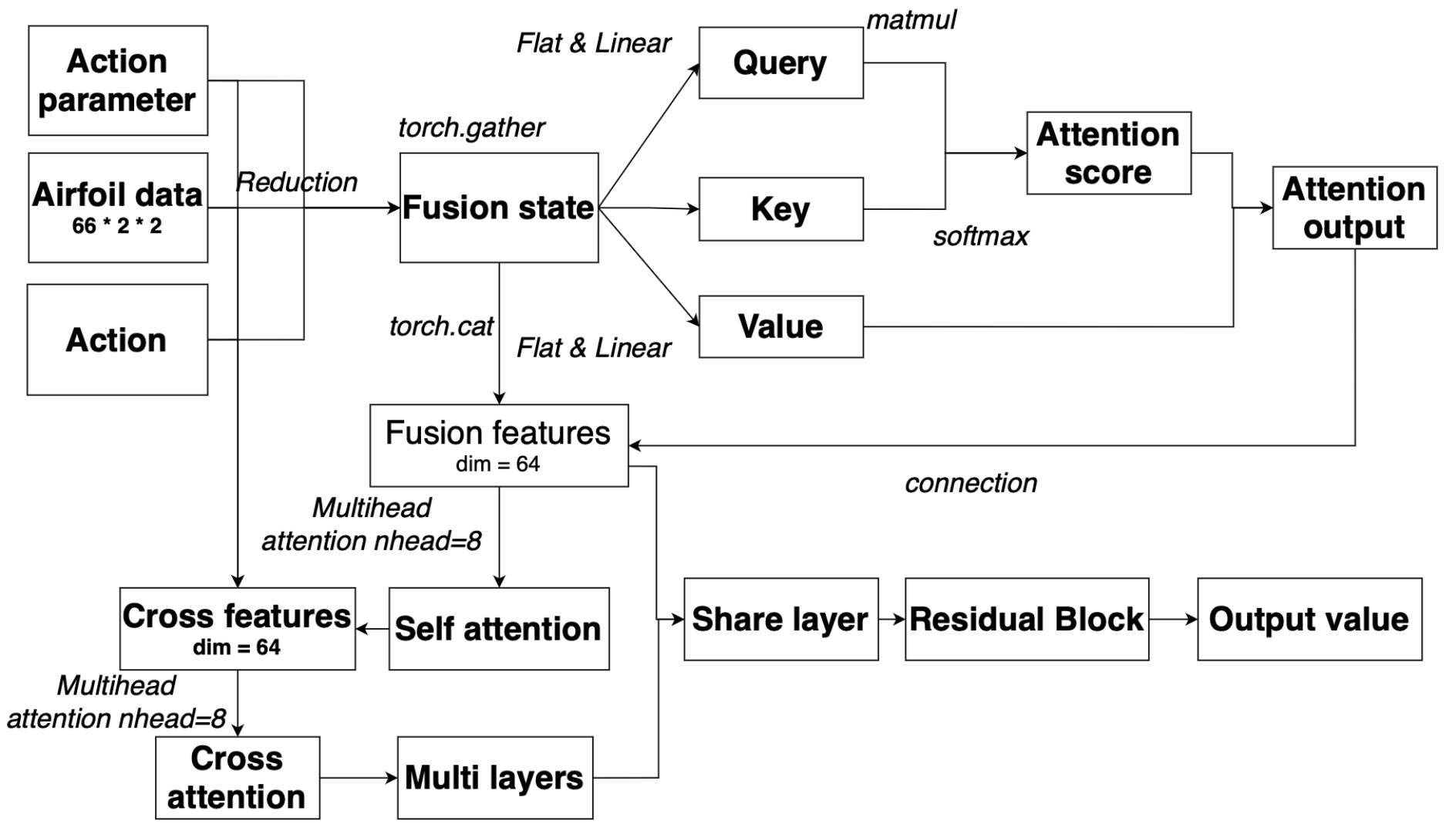}
      \caption{Critic network architecture}  
      \label{tab:critic_network_architecture}
      \end{figure}

The dimension is reduced efficiently with the Bézier curve representation. Then the information should be processed through the neural network to obtain the action that maximizes the cumulative reward from the current state. The output includes action probability for different types of action and corresponding parameters that action possesses. Action that is chosen to be performed originates from the $softmax$ layer while the corresponding parameters are based on the restriction of shape design. For fear of mesh tangling, an activation function needs to be designed that maps the parameters to reasonable intervals. However, the traditional activation functions like $tanh$ and $sigmoid$ suffer from the gradient disappear phenomenon. Then the $softsign$ type activation function is adopted. Techniques like residual connections\cite{he2016deep}, the SELU activation function\cite{klambauer2017self}, and real-time monitoring via TensorBoard\cite{abadi2016tensorflow} were employed to mitigate and monitor potential gradient vanishing problems.

The critic neural network, as illustrated in \cref{tab:critic_network_architecture}, is designed to evaluate the potential feedback for various actions acted on different states based on the experience accumulated in the replay buffer. The reward mechanism is complicated since the data originate from a CFD solver, not even to mention the shape deformation. To approximate such a process, the interaction between the shape, actions and corresponding parameters should be modeled efficiently. We consider adopting the self-attention and cross-attention mechanism together to investigate the inherent relation. The self-attention mechanism can detect those features which take a significant influence on the target functional. Meanwhile, the cross-attention mechanism is adopted to approximate the shape deformation process. It tries to understand how the action with corresponding parameters influenced the shape deformation. The state, action, and action parameters are mapped to the hidden layer with the same dimension to ensure the cross-attention layer can be modeled.

The loss function is an important component which decides the update direction of parameters in neural networks. For different modules, the loss function is defined according to the features. To choose the action that maximizes the reward based on specific states, the agent should evaluate the potential feedback correctly at first. It indicates that the critical neural network should be updated more frequently than the actor neural network. The long-term reward loss is defined as
\begin{equation}
   Critic\_loss = {R}_{t+1}(S_t, A_t| \theta^c) + \gamma * Q(S_{t+1}|\theta^a, \theta^c) - Q(S_{t}|\theta^a, \theta^c)
\end{equation}
is adopted to minimize the error of long-term reward based on the real reward of the environment. The selection of the loss function for the actor network is designed to effectively guide the policy towards actions that maximize expected rewards. This is achieved through a gradient ascent approach on the expected return. The maximal reward loss is typically defined as the negative of the Q-value estimated by the critic network for the current state and the action proposed by the actor network, i.e.,
 \begin{equation}
     Actor\_loss = - Q(S_{t}|\theta^a, \theta^c).
 \end{equation}
 
With the actor and critic neural networks, the agent interacts with the environment and tries to accumulate experience to perform appropriate actions based on different states. The learning performance is highly related to data originating from the interaction process. In order to enhance efficiency, the exploration should incorporate sufficient states. Otherwise, the exploitation may get trapped into overfitting issues. While the shape optimal design of airfoil is a complicated process and not easy to be trained for solving specific target functional, we ameliorate the exploration and exploitation process as follows.
\subsection{Exploration and exploitation}
Exploration and exploitation are fundamental components of the interaction process.
Balancing these two components is vital in reinforcement learning. Excessive exploration is time-consuming. On the other hand, too much exploitation might cause missing out on better strategic options. In this framework, we adopt the $\epsilon$-greedy strategy. This is a popular technique that chooses an action at random with probability $\epsilon$ and selects the action that is predicted to have the highest reward based on its current knowledge with probability $1-\epsilon$. Besides, in our framework, we integrate Ornstein-Uhlenbeck (OU) noise as a complementary tool to the $\epsilon$-greedy strategy, specifically targeting the exploration aspect of reinforcement learning in continuous action spaces. While the $\epsilon$-greedy strategy efficiently balances exploration and exploitation by randomly selecting actions with a probability $\epsilon$, OU noise brings more possibility to the exploration process.

Even though the optimization mechanism can be well implemented with the exploration and exploitation process, the learning effect is not ideal usually. The difficulties originate from the complex mechanisms involved in airfoil design. Firstly, it is not easy to identify the correct reward direction in random exploration. Excessive penalization can lead the neural network parameters into an unfavorable bottleneck, where the agent opts for minimal action to avoid penalties, eventually leading to a suboptimal local optimum.
Secondly, airfoil design inherently demands a geometry with a certain level of smoothness. If the shape changes drastically in a single step, even a correct design trend may result in unfavorable reward feedback due to geometric roughness. This also indicates that the shape optimal design needs to be implemented within an elaborate parameter representation environment. 

To address the above issues, we introduced several improvement strategies, i). if a certain action results in a penalty, the agent attempts to explore the opposite direction starting from the pre-deformation geometry. Both exploration results and parameters are stored in the replay buffer, and the action leading to a better outcome is chosen for the next step. This mechanism significantly reduces aimless exploration in airfoil design, preventing the model from extensively exploring in unproductive directions. This is also inspired by the idea of a twin network to prevent oscillation during the simulation and it turned out to be very effective for reinforcement learning. Similarly, if the current action brings a perfect reward, the agent will be encouraged to take a similar action with a tiny noise. As the environment is not straightforward to obtain positive meaningful rewards, such a mechanism can help the agent to accumulate useful experience quickly. ii). the single step movement of the sample points is constrained, so that there may not occur irregular geometry. Besides, an early stop mechanism is constructed when the accumulated reward decreases to the minimum threshold. This can be explained as the agent failing for the current design process when the geometry gets trapped into a very terrible shape and may not recover easily. Then the environment should be reset for the next turn. iii). the training process is separated into different periods. A small batch size is adopted for the initial stage. When the replay buffer accumulates enough data, a larger batch size number will be updated.
Besides, the data are stored in the replay buffer with a 4-tuple form transition $(state, action, action\_parameters, next\_state)$. Every time the gradient should be calculated from the 4-tuple with a number of batch sizes. The components include three parts, the recent transitions, the transitions with the best reward and randomly chosen transitions from the whole memory. The recent transitions represent the behavior under the current parameters of neural networks while the best transitions can help the agent to repeat the beneficial actions and the random transitions can make the gradient calculation more stable. With such efficient exploration and exploitation mechanisms, we step further to fulfill the reinforcement learning framework.

\subsection{Reward function}

In the context of shape optimal design, a well-designed reward function can guide the reinforcement learning algorithm toward more effective and meaningful solutions, ensuring the learned policies are not only efficient but also aligned with the intended outcomes. Even though techniques introduced during the reinforcement learning process, such as discount factors, might prevent strict equivalence in real computational scenarios, maintaining this equivalence as closely as possible facilitates the learning process. Besides, the mechanism is of vital importance to guarantee a precise value of reward where the agent can learn effective policy.

In \cite{viquerat2021direct,chen2023study,dussauge2023reinforcement}, different reward functions are adopted to meet the requirements of different optimization tasks. It should be noted that the reward function should be carefully designed to support the learning process.
A misalignment between the learning objectives and the actual goals, while technically successful within the learning framework, can lead to solutions that cannot fulfill the intended purpose. The resulting policies might only be effective in specific environment and perform poorly when faced with new scenarios. To overcome above issues, we design reward function and analyze the equivalence of the optimization result.

Let us consider the objective of generating a shape with the global minimum drag value. Denote drag as $D(\{sp^t_i\}_{i=0}^{L+1})$, where $\{sp^t_i\}_{i=0}^{L+1}$ represents the sample points of airfoil at time step $t$ and $L+2$ is the number of sample points as previously defined.  $D(\{sp^*_i\}_{i=0}^{L+1})$ represents one of the global minimum drag values of an airfoil with corresponding parameters. During the optimization process, a sample points sequence will be obtained, $\{sp^{0}_i\}_{i=0}^{L+1}, \{sp^{1}_i\}_{i=0}^{L+1},\cdots \{sp^{t}_i\}_{i=0}^{L+1},\cdots$. While seeking for the objective shape of airfoil, the sequence is expected to converge to a sample points vector equivalently.

\begin{theorem}[Equivalence of reward function] 
\label{theoremRF}
Suppose that $D(\{sp^*_i\}_{i=0}^{L+1})$ represents the global minimum drag value.
If the reinforcement learning process aims to maximize the cumulative reward $Q_t=\sum_{s=t}^{+\infty} R_{s}$ for each step $t$ and generate a sample points sequence with the limit $\{sp^{**}_i\}_{i=0}^{L+1}$, then the result is one of the global minimum of the optimization objective, i.e., $D(\{sp^*_i\}_{i=0}^{L+1}) = D(\{sp^{**}_i\}_{i=0}^{L+1})$.
\end{theorem}

\begin{proof}
    As $D(\{sp^*_i\}_{i=0}^{L+1})$ represents the global minimum drag value of airfoil, then
    \begin{equation}
    \label{defReward}
\begin{array}{cc}
    &D(\{sp^*_i\}_{i=0}^{L+1}) \leq D(\{sp^t_i\}_{i=0}^{L+1})\qquad\forall sp^t_i\in \mathbb{R}^2, \forall t \in \mathbb{N} \\
    &\{sp^*_i\}_{i=0}^{L+1} \equiv \{sp^{t_0}_i\}_{i=0}^{L+1} \quad \text{if~} 
    D(\{sp^{t_0}_i\}_{i=0}^{L+1}) \leq D(\{sp^*_i\}_{i=0}^{L+1}).
    \end{array}
\end{equation}

 Since $\{sp^{**}_i\}_{i=0}^{L+1}$ denotes the sample points of the airfoil that yield the maximal reward of the learning process, then
 \begin{equation}
    D_0 +(-D(\{sp^{**}_i\}_{i=0}^{L+1}))  = \lim\limits_{t\to+\infty}\sum_{s=1}^t R_{s} \geq D_0 + (-D(\{sp^t_i\}_{i=0}^{L+1}))\qquad\forall sp^t_i\in \mathbb{R}^2, \forall t \in\mathbb{N}. 
\end{equation}
More specifically,
\begin{equation}
    D_0 +(-D(\{sp^{**}_i\}_{i=0}^{L+1})) \geq D_0 + (-D(\{sp^*_i\}_{i=0}^{L+1})).
\end{equation}
The equation above can be converted to:
\begin{equation}
     D(\{sp^{**}_i\}_{i=0}^{L+1}) \leq D(\{sp^*_i\}_{i=0}^{L+1}).
\end{equation}
From \cref{defReward}, $\{sp^{**}_i\}_{i=0}^{L+1}$ is equivalent to $\{sp^{*}_i\}_{i=0}^{L+1}$.
\end{proof}
 
The process above demonstrates the equivalence of the reinforcement optimization objective to the primary task. A similar analysis, such as maximizing the lift-drag ratio, can also be derived in this manner.

To further refine the optimization strategy, regularization terms can be added to the reward function. These terms are designed to encourage exploration. By incorporating exploration-driven regularization, the reward function can be represented as:
\begin{equation}
\tilde{R}_{t} = (D_{t-1} -D_t) + \lambda^{\epsilon}_{t} (D_{0}-D_{t}),
\end{equation}
where $\lambda^{\epsilon}_t$ is a weighting coefficient that balances the original objective with the exploration incentive. When the current shape yields an improvement in the objective function value, the agent gets a reward. Conversely, a penalty is incurred for generating a geometry that does not meet the expected criteria. 
The parameter $\lambda^{\epsilon}_t$ shall be decreased with the update of noise. 
In order to guarantee the equivalence above, the coefficient $\lambda^{\epsilon}_t$ is set as an exponential decline function.
\begin{theorem}[Equivalence of generalized reward function] 
If the learning aims to maximize the cumulative reward $\tilde{Q}_t=\sum_{s=t}^{+\infty} \tilde{R}_{s}$ and generates a sample points sequence with limit $\{sp^{***}_i\}_{i=0}^{L+1}$, the result is one of the global minimum of the optimization objective, i.e., $D(\{sp^*_i\}_{i=0}^{L+1}) = D(\{sp^{***}_i\}_{i=0}^{L+1})$.
\end{theorem}
\begin{proof}
    Suppose that the learning process generate a sequence $(D^{**}_t)_{t\in\mathbb{N}}$ which has a limit $D(\{sp^{**}_i\}_{i=0}^{L+1})$, if the objective is to maximize the cumulative reward $\sum_{s=t}^{+\infty} {R}_{s}.$
Similarly,  a sequence $(D^{***}_t)_{t\in\mathbb{N}}$ is generated with the objective to maximize the cumulative reward $\sum_{s=t}^{+\infty} \tilde{R}_{s}.$
 Firstly, the cumulative reward is finite since the coefficient is an exponential decline function, i.e.,
    \begin{equation}
        \sum_{s=1}^{+\infty} \tilde{R}_{s} = \sum_{s=1}^{+\infty}  {R}_{s} + \sum_{t=1}^{+\infty}  \lambda^{\epsilon}_t (D_0^{***}-D_t^{***}) \leq D_0^{***} + (D_0^{***}) \sum_{t=1}^{+\infty}  \lambda^{\epsilon}_t < +\infty.
    \end{equation}
If $D(\{sp^{***}_i\}_{i=0}^{L+1})\leq D(\{sp^{**}_i\}_{i=0}^{L+1}),$ from the definition and \cref{theoremRF}, the equivalence of reward function held naturally. To the contrary, suppose $D(\{sp^{***}_i\}_{i=0}^{L+1}) > D(\{sp^{**}_i\}_{i=0}^{L+1})$ and denote $\delta = ( D(\{sp^{***}_i\}_{i=0}^{L+1})- D(\{sp^{**}_i\}_{i=0}^{L+1}))/2$. Then, $\exists N_{\delta}$, s.t., $\forall t > N_{\delta}$,  $D_t^{***} - D_t^{**} > \delta /2$. While the sequence $(D^{***}_t)_{t\in\mathbb{N}}$ and $(D^{**}_t)_{t\in\mathbb{N}}$ are generated with maximizing $Q_{t}$ and $\tilde{Q}_t$, specifically, we have
\begin{equation}
\begin{aligned}
\label{inequal}
  &+\infty > \tilde{Q}_{N_{\delta}} = \sum_{t=N_{\delta}}^{+\infty} \tilde{R}_{s}  = \sum_{t=N_{\delta}}^{+\infty} \left((D_{t-1}^{***}-D_t^{***})   +  \lambda^{\epsilon}_t (D_0^{***}-D_t^{***})\right) \\
    &\geq \sum_{t=N_{\delta}}^{+\infty} \left((D_{t-1}^{**}-D_t^{**})   + \lambda^{\epsilon}_t (D_0^{**}-D_t^{**})\right)
    \geq    \sum_{t=N_{\delta}}^{+\infty} \left((D_{t-1}^{***}-D_t^{***}) +  \lambda^{\epsilon}_t (D_0^{**}-D_t^{**})\right). 
\end{aligned}
\end{equation}
Since the sequences both starts from $D_0$ and share the same coefficient $\lambda^{\epsilon}_t$, 
the following inequality can be derived
\begin{equation}
 \sum_{t=N_{\delta}}^{+\infty}  \lambda^{\epsilon}_t (D_t^{**}-D_t^{***}) \geq 0.
\end{equation}
However, the inequality cannot held since $D_t^{***} - D_t^{**} > \delta /2, ~\forall~t > N_{\delta}$, which implies the equivalence of this generalized reward function is held.
\end{proof}

Thus this generalized reward function will encourage exploration at the initial period while the result can converge to the optimization objective with the iteration getting processed. 
These mechanisms ensure a more robust and effective search for the optimal airfoil shape, facilitating the discovery of innovative designs that minimize drag while navigating the complex landscape of possible configurations.

\subsection{Twin delayed DDPG algorithm}
With the introduced modules, we are going to develop a more robust framework,
TD3. It introduces several key improvements that enhance the stability and performance of learning tasks. The primary advantage of TD3 is its twin-critic network architecture. It involves two separate critic networks to estimate the Q-values, where the smaller one is used for the policy update. This approach effectively reduces the overestimation bias often seen in DDPG, leading to more stable and reliable learning.
Additionally, TD3 introduces policy update delays, where the policy network is updated less frequently than the critic. This delay helps mitigate the risk of the policy network exploiting inaccuracies in the Q-value estimates, further enhancing the stability of the learning process. With such a more stable framework, we constructed Algorithm \ref{td3alg}.
\begin{algorithm}
\caption{Twin Delayed Deep Deterministic Policy Gradient}
\label{td3alg}
\begin{algorithmic}[1]
\REQUIRE Initial airfoil shape data, actor parameters $\theta^a$, critic parameters $\theta^{c_1}$ and twin critic parameters $\theta^{c_2}$, target networks  $\theta^a_{target}$,  $\theta_{target}^{c_1}$ and  $\theta_{target}^{c_1}$, ReplayBuffer
\ENSURE Trained Twin Delayed DDPG agent
\STATE Set the target parameters $\theta^a_{target} \leftarrow \theta^a, \theta^{c_1}_{target} \leftarrow \theta^{c_1}, \theta^{c_2}_{target} \leftarrow \theta^{c_2}$
\FOR{each episode in max\_episodes}
    \IF{Exploration or reach the minimum threshold}
        \STATE Reset the environment and get the initial state
    \ELSE
        \FOR{step in episode}
            \IF{Exploration}
                \STATE The agent produces a random action for execution
            \ELSE
                \STATE The agent selects an action based on $\theta^a_{target}$
                \STATE Noises are added to the action and corresponding parameters
                \STATE The agent executes the action based on the current state
                \STATE Update the cache state
            \ENDIF
            \STATE Observe the new state and calculate the reward based on AFVM4CFD
            \IF{reward $<$ min\_reward}
                \STATE Reload the cache state
                \STATE Execute the opposite action with corresponding parameters
            \ENDIF
            \STATE Save the tuple (state, action, reward, next\_state) to the ReplayBuffer
            \IF{ReplayBuffer has enough experiences}
                \STATE Sample experiences (Recent, Best, Random) from ReplayBuffer
                \STATE Compute the target actions based on $\theta^a_{target}$
                \STATE Compute the target Q-value $\min(Q_{\theta_{target}^{c_1}}, Q_{\theta_{target}^{c_2}})$
                \STATE Update the critic networks $\theta^{c_1}$, $\theta^{c_2}$
                \STATE Soft update critic networks $\theta_{target}^{c_i} \leftarrow \tau \theta_{target}^{c_i} + (1-\tau) \theta^{c_i}$ for $i = 1, 2$
                \IF{$step \equiv 0 \pmod{delay\_step}$}
                    \STATE Update the actor network $\theta^{a}$
                    \STATE Soft update actor networks $\theta_{target}^{a} \leftarrow \tau  \theta_{target}^{a} + (1-\tau)\theta^{a}$
                \ENDIF
            \ENDIF
            \STATE Update the state
        \ENDFOR
    \ENDIF
\ENDFOR
\end{algorithmic}
\end{algorithm}

\section{Numerical Result}
\label{sec:result}
\subsection{Drag minimization}
The initial example focus on minimizing drag calculations. The study utilizes a NACA0012 airfoil, positioned within a circular domain with a radius of $35$, under a boundary condition of mirror reflection. The operational conditions are set with Mach number $0.85$ and attack angle $0^\circ$.

 \begin{figure}[!h]\centering
        \includegraphics[width=0.48\textwidth]{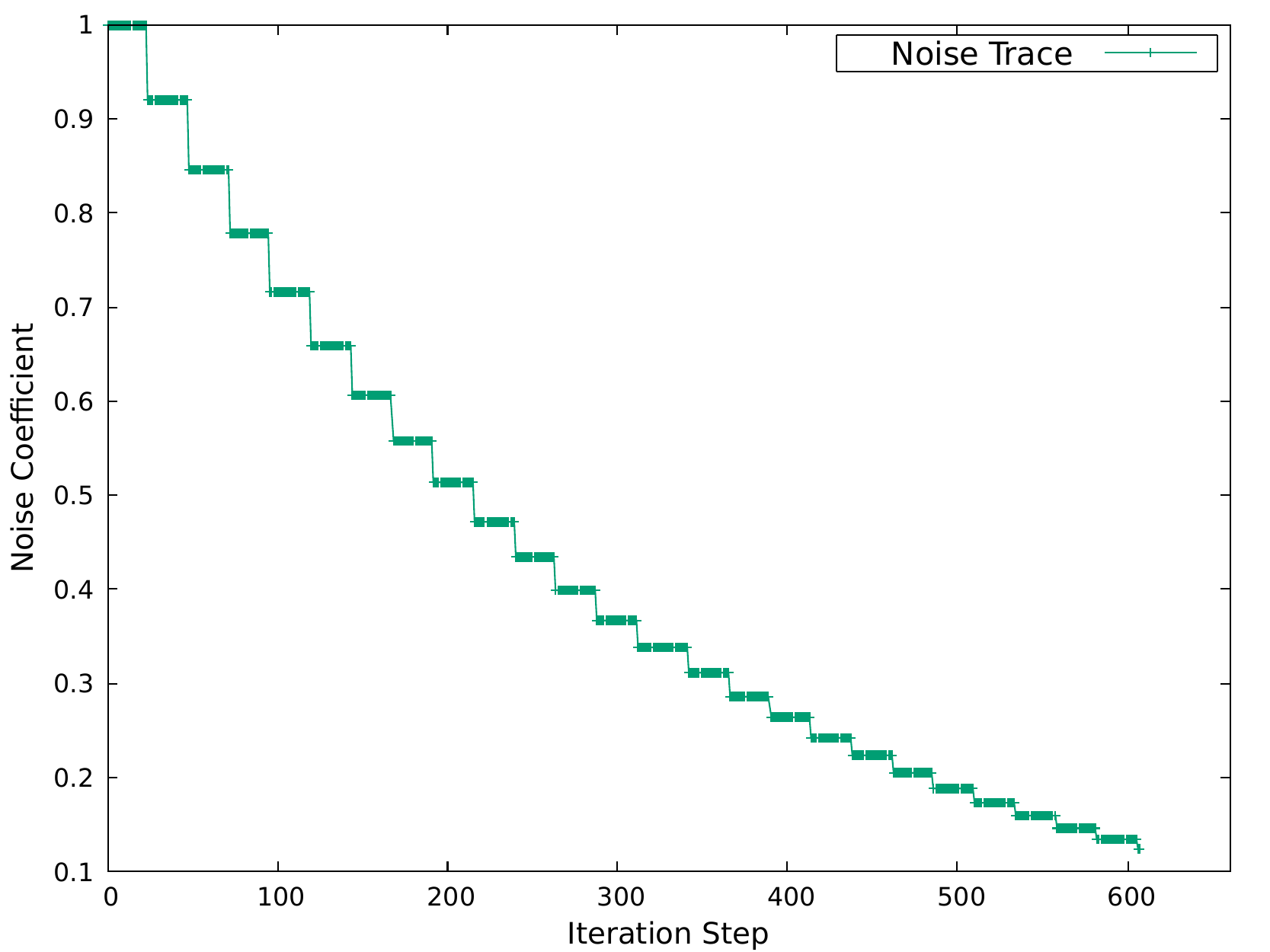}
        \includegraphics[width=0.48\textwidth]{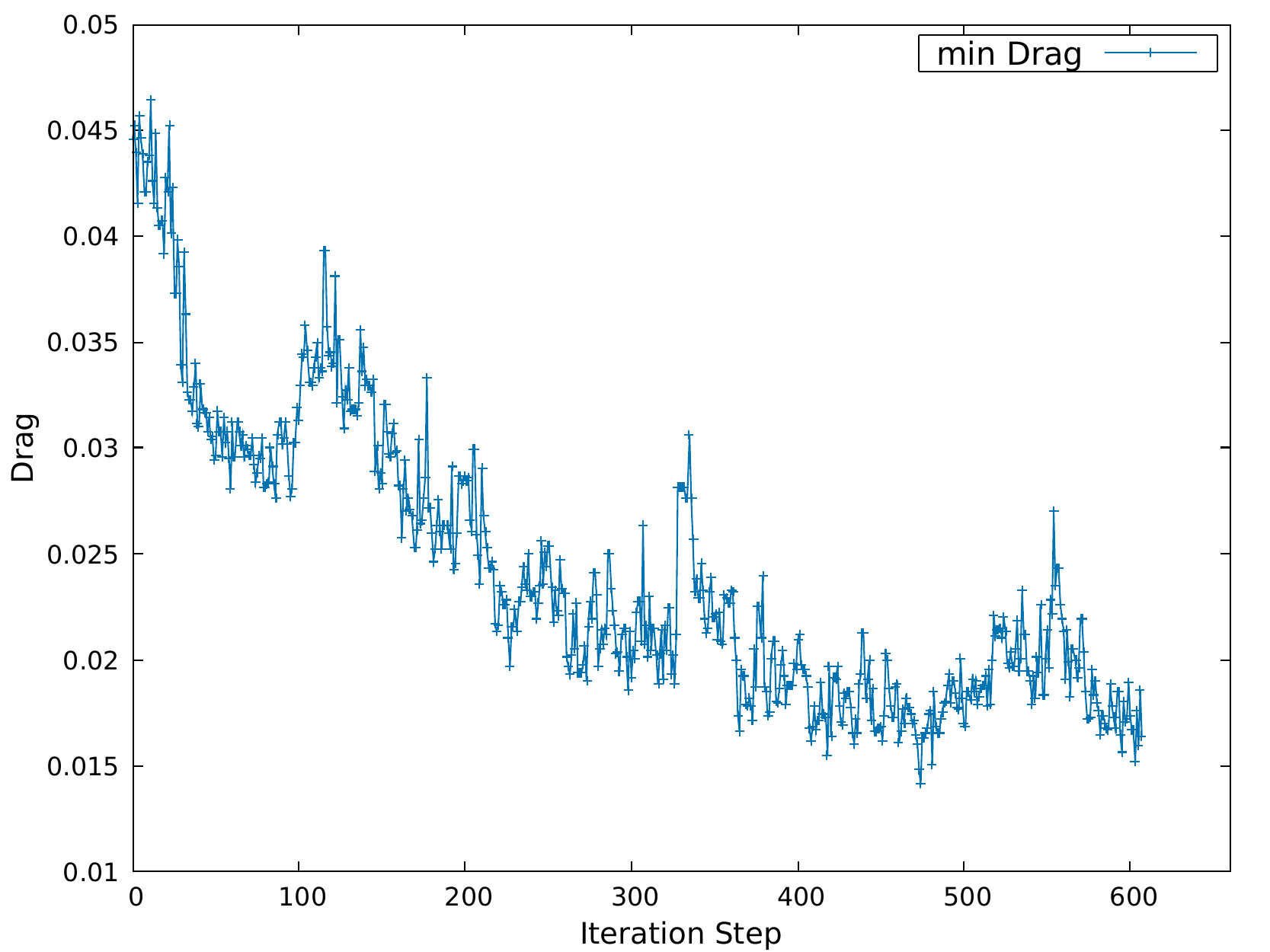}
       \caption{Record curves during the optimization. Left: noise coefficient. Right: drag value.}
      \label{noiseTrace}
      \end{figure}
      
To ensure the integrity of the airfoil's structure during training, a constraint is imposed, mandating a minimum thickness of 0.1 for sections where $x<0.1$ and $x>0.7$. Should a potential action result in an airfoil design adhering to this thickness constraint, it is deemed physically infeasible. Consequently, this action will not be executed, leading to a negative reward.

      As shown in Algorithm 3, the agent embarks on accumulating experience through random actions performed four times, each starting from the initial configuration across $64$ time steps. The initial recorded drag value is $0.0452106$.
\begin{figure}[h]\centering
      \frame{\includegraphics[width=0.32\textwidth]{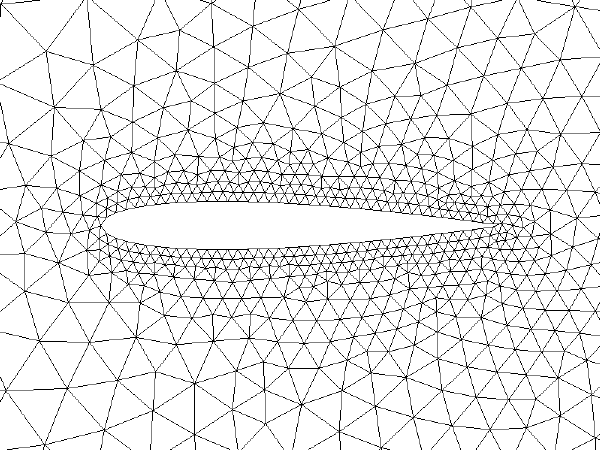}}
        \frame{\includegraphics[width=0.32\textwidth]{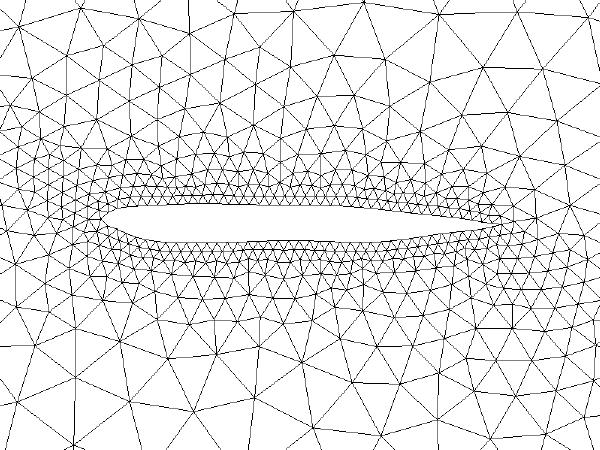}}
        \frame{\includegraphics[width=0.32\textwidth]{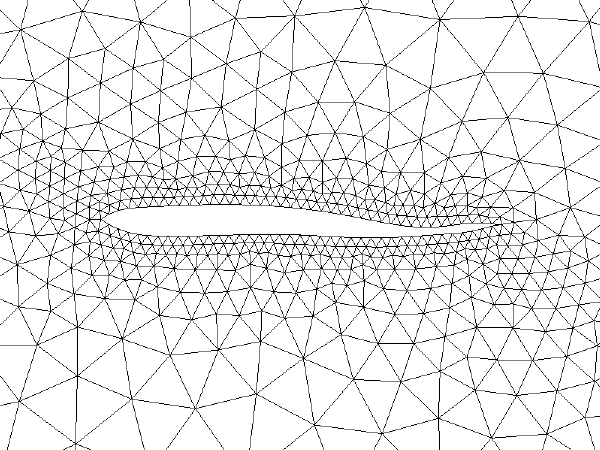}}
       \caption{Right: initial airfoil, NACA0012. Middle: the airfoil generated from the 15th episode. Right: the result of the optimized airfoil.}
      \label{minDrag}
      \end{figure}
Training starts with a moderate batch size of $256$, spanning a total of $75$ epochs, with each epoch comprising 8 time steps. The noise coefficient gets updated every three epochs, following an exponential decline. As illustrated in Figure \ref{noiseTrace}, the drag value exhibits a decrement trend as iterations advance. Even though the record curve is not smooth as shown in the record line, it still evolved towards the optimization objective. The non-smooth phenomenon can be explained by the fact of noise.
Even if the potential action from the actor neural network generates an effective action with corresponding parameters, the noise may distract the agent from the potential actions. However, such a noise module is indispensable to help the agent explore the possibility of different states.

Throughout the optimization process, the agent explores various airfoil shapes. For instance, Figure \ref{minDrag} exhibits a configuration from the $15^{th}$ epoch, with a drag value of $0.031537$. The drag value's trajectory halts around $0.02$, constrained by the minimum thickness. Encounters with this restriction trigger a negative reward, arousing exploration toward alternate states. Notably, as the exploration progresses, particularly in the concluding epoch, the airfoil's tail approaches the thickness restriction.

The interconnection between drag optimization and shock waves is a fundamental aspect of airfoil design, where minimizing the impact of shock waves is pivotal for reducing wave drag. This concept has been supported and observed across various studies within the field \cite{mazaheri2015optimization}. From our experiment, the NACA 0012 airfoil initially exhibited strong shock waves on both the upper and lower surfaces, typically associated with substantial pressure discontinuities near these regions, significantly contributing to increased drag. Through the optimization process, notable changes in the structure of these shock waves were observed. Specifically, the optimization resulted in the attenuation and positional shift of the shock wave on the airfoil's upper surface, along with the complete dissipation of the lower surface's shock wave.

      \begin{figure}[h]\centering
        \frame{\includegraphics[width=0.48\textwidth]{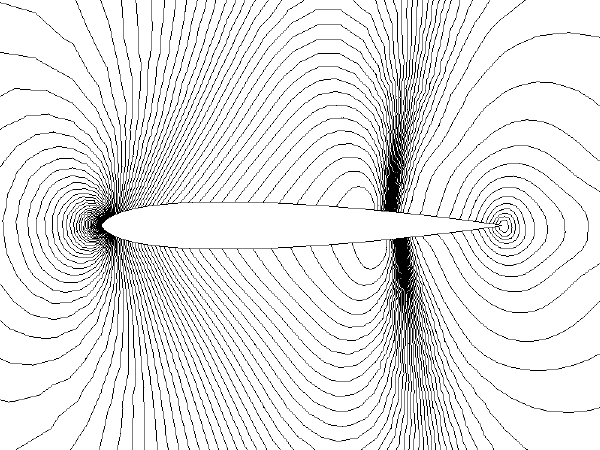}}
        \frame{\includegraphics[width=0.48\textwidth]{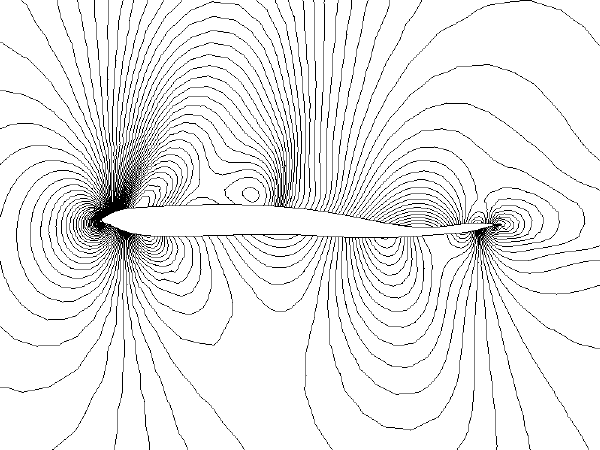}}
         \caption{Mach isoline. Left: the initial airfoil. Right: optimized airfoil.}
      \label{Machline}
      \end{figure}
\subsection{Lift-drag ratio maximization}

      \begin{figure}[h]\centering
        \includegraphics[width=0.492\textwidth]{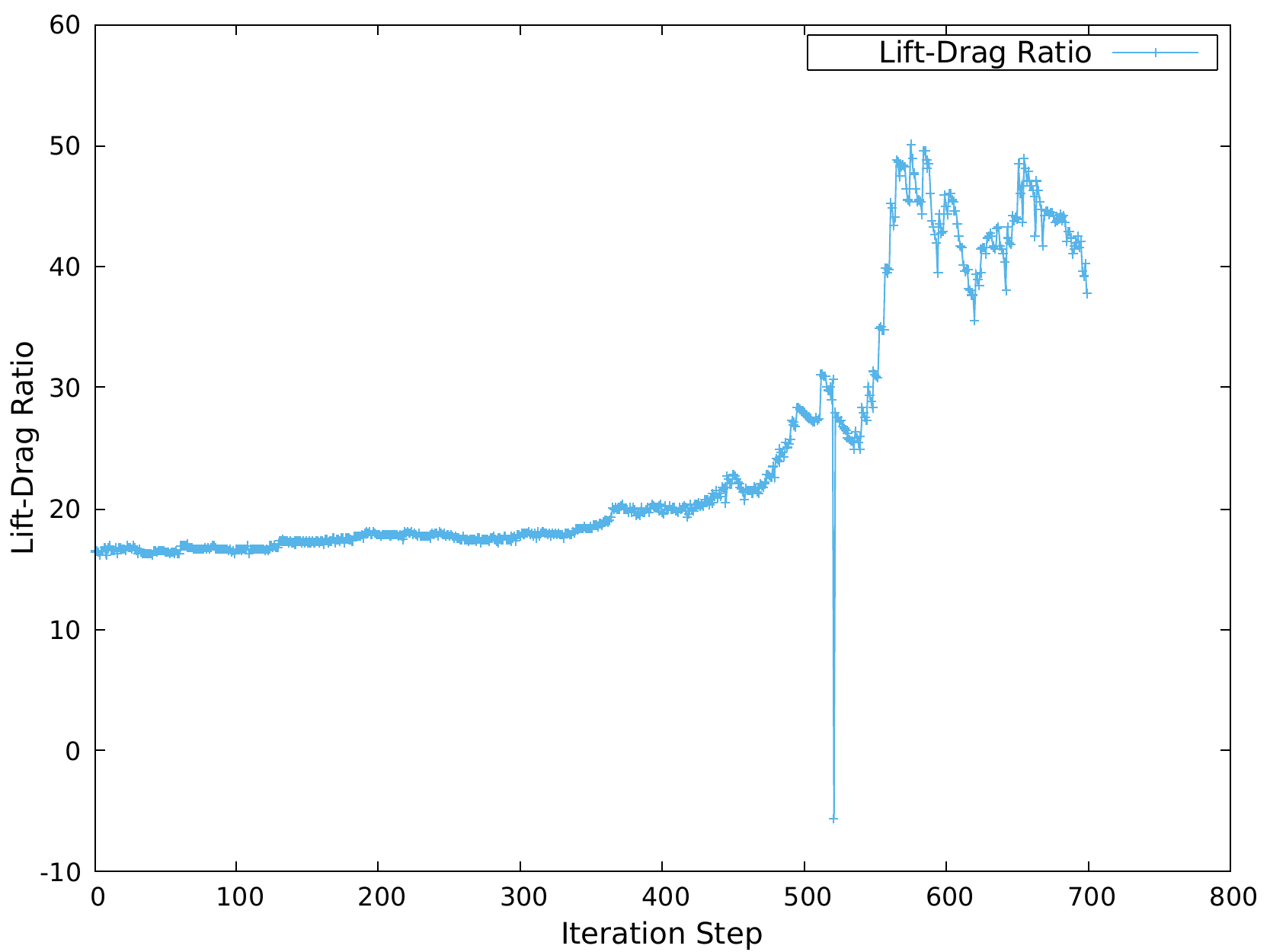}
        \frame{\includegraphics[width=0.48\textwidth]{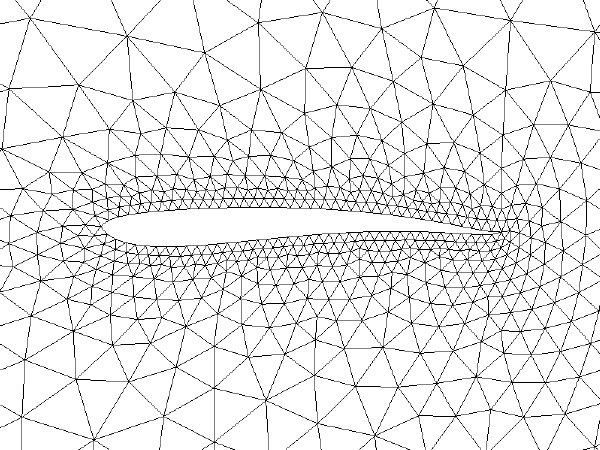}}
         \caption{Left: record curve of lift-drag ratio during the optimization process. Right: mesh around the optimized airfoil.}
     
      \label{ratioRecord}
      \end{figure}
      
    We conducted the optimization task of maximizing the lift-drag ratio of NACA0012 with Mach number $0.8$ and attack angle $1.25^\circ$. The value of the lift-drag ratio gets raised from $16.488$ to $48.916$.
    Optimizing for the lift-drag ratio poses greater challenges than solely focusing on drag due to the singularity of the objective functional. Firstly, the calculation of lift is more complicated than drag. Evidence from our previous experiments and other researchers \cite{DOLEJSI2021178} have shown that the convergence curve of lift during the calculation is oscillatory. This may be caused by the weaker regularity of the dual solution for the lift coefficient. 
    Secondly, the DWR-based mesh adaptation method can not be generalized to the lift-drag ratio easily. Then, the accuracy of the target functional gets influenced a lot. Multi-mesh method\cite{kuang2023mul} can be applied to solve the lift and drag respectively. However, it will bring additional difficulties to the PDE solver.

          \begin{figure}[t]\centering
        \frame{\includegraphics[width=0.48\textwidth]{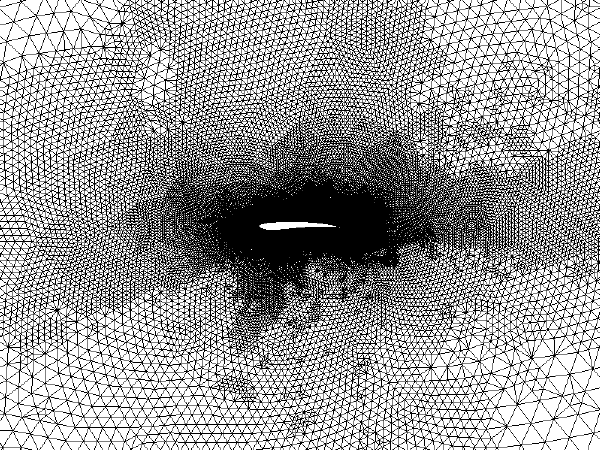}}
        \frame{\includegraphics[width=0.48\textwidth]{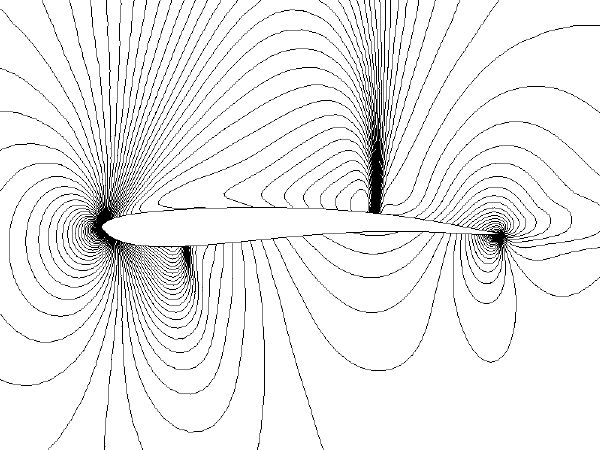}}
          \caption{Left: mesh generated from DWR-based adaptation. Right: Mach isoline.}

      \label{ratioMach}
      \end{figure}
    Currently, we pay more attention to the drag calculation since the denominator demands higher accuracy in the lift-drag ratio calculation. The mesh generated from DWR-based adaptation as shown in Figure \ref{ratioMach} is adopted to calculate the lift-drag ratio. 
    Our experiments, illustrated in Figure \ref{ratioRecord}, highlight  data anomalies issue that can arise when the framework for solving the objective functional is not sufficiently accurate. Such inaccuracies can lead to erroneous signals that affect the policy update process. It indicates that a more task-specific design of the reward functional need to be designed. The optimization process, depicted in Figure \ref{ratioRecord}, shows an extensive period of exploration to identify a viable direction for improvement. Although the optimized shape has a lift-drag ratio nearly three times than the original one, the airfoil shape does not represent the intended airfoil design. The agent should explore more potential shapes to generate the expected geometry, which will consume time much longer than the current experiment. While the current PDE solver has already been a highly efficient one, an increase of magnitude with respect to the wall time is tough. Multi-agent exploration framework should be constructed in the future to enhance this algorithm. Shock waves can be found in Figure \ref{ratioMach}, which may explain the reason that the lift-drag ratio stopped around $45$ during the optimization process. To develop a task-specific design algorithm for the lift-drag calculation, we need to pay attention to the pressure distribution along the geometry as well.

\section{Conclusion}
\label{sec:conclusion}
      In this work, we present a mechanism-driven reinforcement learning framework. It integrates reinforcement learning with goal-oriented PDEs-based systems for shape optimization of airfoil.
      We developed a suite of mesh recovery techniques aimed at improving mesh quality to overcome the challenge brought by high-dimensional states and unstructured meshes.
      It is important to highlight that the reinforcement learning approach was meticulously tailored to align with the intrinsic motivation behind the optimization objectives.
Upon conducting experimental validations of shape optimization, our findings demonstrate that the proposed framework is adept at managing simulations with precisely calculated objective functions. 

However, more complex target functionals have not yet been further explored. Looking forward, we aim to integrate advanced techniques, such as the multi-mesh method for tackling more intricate target functionals and employing multi-agent strategies for enhanced exploration efficiency. It should be noted that multi-agent systems can achieve a deeper integration with traditional optimization algorithms based on adjoint equations, leading to more efficient exploration patterns. Furthermore, the adoption of various acceleration techniques, including the incorporation of additional learning modules, parallelization strategies, and more efficient network architectures, represents a promising direction for future development.

\section*{Acknowledgments}
The authors would like to thank the support from The Science and Technology Development Fund, Macao SAR (No. 0082/2020/A2), National Natural Science Foundation of China (Nos. 11922120, 11871489), MYRG of University of Macau (No. MYRG-GRG2023-00157-FST-UMDF).

\bibliographystyle{unsrt} \bibliography{rl4op}
\end{document}